\documentclass[reqno]{siamart190516}
\usepackage{braket,amsfonts}
\usepackage{amssymb}
\usepackage{mathrsfs}
\usepackage{array}
\usepackage{stmaryrd}
\usepackage{titlesec}
\usepackage{enumerate}
\usepackage{clipboard}
\usepackage{url}
\newclipboard{myclipboard_main}
\usepackage{xr}
\titleformat{\paragraph}
{\normalfont\normalsize\bfseries}{\theparagraph}{1em}{}
\titlespacing*{\paragraph}
{0pt}{3.25ex plus 1ex minus .2ex}{1.5ex plus .2ex}
\usepackage{algorithm}
\usepackage{algpseudocode}
\usepackage{algorithmicx}
\makeatletter
\newenvironment{breakablealgorithm}
{
	\begin{center}
		\refstepcounter{algorithm}
		\hrule height.8pt depth0pt \kern2pt
		\renewcommand{\caption}[2][\relax]{
			{\raggedright\textbf{\ALG@name~\thealgorithm} ##2\par}%
			\ifx\relax##1\relax 
			\addcontentsline{loa}{algorithm}{\protect\numberline{\thealgorithm}##2}%
			\else 
			\addcontentsline{loa}{algorithm}{\protect\numberline{\thealgorithm}##1}%
			\fi
			\kern2pt\hrule\kern2pt
		}
	}{
		\kern2pt\hrule\relax
	\end{center}
}
\makeatother
\usepackage{graphics}
\usepackage{graphicx}
\usepackage{caption,subcaption}
\usepackage{multicol}

\newsiamthm{lem}{Lemma}
\newsiamthm{prop}{Proposition}
\newsiamthm{thm}{Theorem}
\newsiamremark{rmk}{Remark}
\newsiamthm{ass}{Assumption}
\newsiamthm{defn}{Definition}
\newsiamthm{pro}{Problem}

\newcommand{\change}[1]{{\color{black}#1}}
\renewcommand{\theequation}{\thesection.\arabic{equation}}
\makeatletter\@addtoreset{equation}{section} \makeatother

\allowdisplaybreaks

\begin{document}
	\title{Randomized Optimal Stopping Problem in Continuous time and reinforcement learning algorithm}
	\thanks{Y. Dong was supported by the National Natural Science Foundation of China (No. 12071333 \& No. 12101458)}
	\author{Yuchao Dong
		\thanks{School of Mathematical Sciences, Tongji University,  Shanghai 200092, China (\email{ycdong@tongji.edu.cn})}
		}
\maketitle
\begin{abstract}
In this paper, we study the optimal stopping problem in the so-called exploratory framework, in which the agent takes actions randomly conditioning on current state and a regularization term is added to the reward functional.  Such a transformation  reduces the optimal stopping problem to a standard optimal control problem. For the American put option model, we derive the related HJB equation and prove its solvability. Furthermore, we give a convergence rate of policy iteration and compare our solution to the classical American put option problem. Our results indicate a trade-off between the convergence rate and bias in the choice of the temperature constant. Based on the theoretical analysis, a reinforcement learning algorithm is designed and numerical results are demonstrated for several models.
\end{abstract}


\section{Introduction}
Reinforcement learning (RL, for short) is about how software agents choose actions in an environment  to achieve some goals or to maximize rewards. Recently, it became one of the most active and fast developing areas of machine learning, due to its success in playing go \cite{silver2016mastering,silver2017mastering}, achieving human-level performance in video games \cite{mnih2015human}, controlling robotics \cite{deisenroth2013survey}, designing autonomous driving \cite{levine2016end} and so on. Applications of RL in the financial industry such as algorithmic trading and portfolio management also have attracted more attention in recent years, see \cite{nevmyvaka2006reinforcement,hendricks2014reinforcement,moody1998performance,moody2001learning} for instance.

One distinguishing feature of  reinforcement learning is that the model or dynamic for the environment may not be known priorly. The agent is not told what to do, but instead, discover which action yields the best result through the interaction with the environment.  This is a case of "kill two birds with one stone": the agent's actions serve both as a mean to explore (learn) and a way to exploit (optimize). Since the exploration is inherently expensive in terms of resources, time and opportunity, the agent must balance between greedily exploiting what has been learned so far to choose actions that yield near-term higher rewards and continuously exploring the environment to acquire more information for potential long-term benefits. Extensive studies have been carried out to find the best strategies for the trade-off between exploration and exploitation. Most  past works do not include exploration as a part of the optimization objective, but treat exploration separately as an ad-hoc chosen exogenous part, see \cite{sutton2011reinforcement}. On the other hand,  \cite{ziebart2008maximum,nachum2017bridging,fox2015taming} propose and apply a discrete-time entropy-regularized RL formulation that incorporates exploration into the optimization objective as a regularization term with a trade-off weight imposed on the entropy on the exploring strategy. Recently, Wang et al \cite{wang2020reinforcement} proposed and developed a continuous-time entropy regularized relaxed control  framework. They called it an exploratory formulation, where the actions are chosen to be probability measures. \Copy{Minor3}{\change{  To be precise, they considered the following controlled dynamic
$$
dX_t=\int_U b(t,X_t,u) \pi_t(u)du dt+\sqrt{\int_U \sigma^2(t,X_t,u)\pi_t(u)du}dW_t,
$$
where the adapted function-valued process $\pi_t(u)$ is a probability density over $U$. The utility functional to be maximized is defined as 
$$
\mathbb E\left[ \int_0^T \left(\int_U f(t,X_t,u)\pi_t(u)du+\lambda H(\pi_t)\right)dt +g(X_T) \right],
$$
with $H(\pi_t):=-\int_U \pi_t(u)\log\pi_t(u)du$ being the differential entropy and $\lambda$ being a constant. We see that, different from classical formulation for optimal control problems,  the entropy is added to the utility functional so as to encourage exploration.   For LQ control problems  with $U$ being $\mathbb R$,  Wang et al \cite{wang2018exploration} show that the optimal distributions of the exploratory control policies must be Gaussian, thereby providing an interpretation for the Gaussian exploration widely used both in RL algorithm design and in practice, see \cite{nachum2017bridging}. For general cases, Tang et al \cite{tang2021exploratory} give the related HJB equation and further studied its properties. }}

In this paper, one motivation is to extend the exploratory framework to the optimal stopping problem. The classical decision problem of optimal stopping has found an amazing range of applications from finance, statistics, marketing, phase transitions to engineering. In many cases, the state space is decomposed into stopping region and continue region and the optimal stopping strategy is to stop whence the state hits the stopping region. These regions are separated by  graph of a function which is also called stopping boundary. Theoretical analysis indicates that this boundary is closely related to the value function of the problem.  In low-dimensional cases, one can accurately compute the boundary by methods based on PDEs and dynamic programming. For high-dimensional cases, one faces the curse of dimensionality and people resort to machine learning techniques. 

\change{ Becker et al \cite{becker2019deep} develop a machine learning method for solving optimal stopping problem. They use neural networks to approximate the stopping time and find the optimal parameters by gradient descent. Later, Felizardo and Matsumoto \cite{felizardo2022solving} improve this method by using convolution neural networks instead of fully-connected ones. Recently, Reppen et al \cite{reppen2022neural} propose another method for solving optimal stopping problem. Different from previous works, the authors use neural networks to approximate the stopping boundary. As indicated in previous paragraph, the optimal strategy is to stop whence the state hits the stopping region. However, from a learning perspective, it is not possible to apply gradient descent directly over these 'stop-or-go' strategies as the lack of smoothness. Thus, the authors use  a soft stopping rule or fuzzy stopping boundary during the learning.  In some sense, it means that the agent stops with a probability related to the distance to the boundary, see \cite{reppen2022neural} for details. Then, one natural question is that what is the effect of such fuzziness on the learning procedure and obtained strategy. In this paper, we will give an insight on this problem that it improves convergence rate and induces bias at the same time. Thus, people has to strike a balance in the choice of randomized policies. }

\Copy{Minor1}{\change{Besides previous works, there are also papers that solving optimal stopping problems arising from different areas with RL methods. These include radiotherapy treatments \cite{ajdari2019towards}, feature selection \cite{liyanage2019automating}, and neural network training \cite{dai2019bayesian}. Fathan and Delage \cite{fathan2021deep} present for the first time a comprehensive empirical evaluation of the quality of optimal stopping policies identified by three state of the art deep RL algorithms: double deep Q-learning, categorical distributional RL, and Implicit Quantile Networks. All of these works consider the discrete time framework and the entropy is not included in the utility functional. In this paper, we will extend the exploratory framework proposed by \cite{wang2018exploration} to optimal stopping problems. We shall see that it is not a trivial extension. }}

\Copy{r2c3}{\change{In mathematics, the value function $u^*$ of an optimal stopping problem is related to the variational inequalities. For example, for American put option, it satisfies the following PDE
\begin{equation*}
	\min\{-\partial_t u^*-\mathcal L_x u^*,  u^*(x,t)-h(x) \}=0,u(x,T)=h(x).
\end{equation*}      
This is one kind of free boundary problems, in which the boundary of the domain where the PDE is defined is part of the solution and needs to be determined as well. For solving high-dimensional PDEs, E and his co-authors \cite{becker2019deep,han2017deep} develope a scheme  for solving a class of quasi-linear PDEs which can be represented by forward-backward stochastic differential equations. Almost at the same time, Raissi et al \cite{raissi2017physics} introduce physics informed neural networks (PINN). They estimate deep neural network models which merge data observations. Their approach solves PDEs arising from physics models in one and two spatial dimensions using deep neural networks. After that, PINN has been further studied and applied to many other problems, see the survey paper \cite{karniadakis2021physics} for more details.  For free boundary problems, Sirignano and Spiliopoulos  \cite{sirignano2018dgm} propose deep Galerkin method that can be used to solve them. Wang and Perdikaris \cite{wang2021deep} also propose a multi-network model based on PINN to tackle a general class of forward and inverse free boundary problems called Stefan problems. By approximating the value function of the optimal stopping problem, we also give an approach to solve free boundary problems. Different from previous methods, ours is based on RL method, which can be model-free, i.e. one does not need the knowledge of coefficients.   

}}
\Copy{r1Contribution}{
\change{Our major contributions are summarized as follows.
\begin{enumerate}
	\item Mathematically, we extend the continuous time exploratory framework proposed by \cite{wang2018exploration} to optimal stopping problems. We introduce an additional state, which represents the surviving probability in some sense. This formulation is similar to that in \cite{touzi2002continuous,laraki2005value}. Different from these works, an regularizer is included in the target functional. Instead of entropy, we introduce another regularizer which makes the problem mathematical tractable. Using the theory of PDEs, we are able to quantify the effect of randomization. Our theoretical result indicates a trade-off between the convergence rate and bias in the choice of the temperature constant, which is also confirmed by numerical experiment.  Although we only consider American option problems, the method of our proof  can be extended to more general models and problems, including Dynkin games and optimal switching. This will be one topic of future research.   
	\item  Our formulation reduces  optimal stopping problems to  standard optimal control problems. With such a reduction, one can adopt RL methods to learn the optimal execution policy for optimal stopping problems. In this paper, we design an offline RL algorithm to learn the optimal execution strategy without model information.  Although the theoretical analysis is under American put option model, the algorithm is applicable for general settings. Our numerical results compare well to that obtained in \cite{becker2019deep,reppen2022neural}. At last, our theoretical results quantifies the effect for the choice of $\lambda$. Combining the sample efficiency analysis, we hope that this will give us an guidance on how to choose proper temperature constant.     
\end{enumerate}
}
}
The paper is organized as follows. In Section \ref{sect_formulation}, we introduce the exploratory formulation. Section \ref{sect_analy} contains the related analytical results on policy iteration and comparison with classical optimal stopping problem. We design a reinforcement learning algorithm and demonstrate some numerical results in Section \ref{sect_rl}. All the proofs are put in Section \ref{sect-proof}.  
\section{Problem Formulation}\label{sect_formulation}
For theoretical analysis, we mainly focus on the American put option problem with one risky asset. Under risk neutral measure, the stock price $S_t$ satisfies the following
$$
dS_t=S_t\left[r dt+\sigma dW_t \right],
$$
with $r$ and $\sigma$ representing risky-free rate and volatility of the stock respectively. Denote by $\mathbb F=\{ \mathcal F_t\}_{0 \le t \le T}$ the filtration generated by Brownian motion $W$. The price of an American put option is related to an optimal stopping problem, i.e. to choose a stopping time $\tau^*$ such that
$$
\mathbb E\left[ e^{-r\tau^*}(K-S_{\tau^*})^+\right]=\sup_{\tau \in \mathcal T[0,T]} \mathbb E\left[e^{-r\tau} (K-S_{\tau})^+\right],
$$
where $\mathcal T[0,T]$ is the totality of all the stopping time  less than $T$. It is well-known that, for such kind problem, the state space can be decomposed into two regions: continue and stopping region, and the optimal stopping strategy is to execute the option whenever state process hits stopping region. 

{\bf Notations} For any region $\Omega$, $C(\Omega)$ represents the set of continuous functions defined on $\Omega$. $C_b(\Omega)$ is the totality of bounded continuous functions equipped with the supreme norm $\|\cdot\|_{\infty}$. $C^{2+\alpha,1+\alpha/2}(\mathbb R \times [0,T))$ presents  a subset of $C(\mathbb R \times [0,T))$ with functions that have $\alpha$-H\"older continuous derivatives up to second order in space and $\alpha/2$-H\"older continuous derivatives up to first order in time. \change{  In this paper, we use $C$ to represent a constant, but may be  different from line to line. In general, it depends on the coefficients $r,\sigma$ and $K$ of the model. In the proof of theorems, we use the notation $C(\cdot)$ to indicate its dependence on other quantities. For example, $C(\varepsilon)$ means  that the constant depends also on $\varepsilon$.}  

\subsection{A Discrete Time Toy Model}\Copy{legtimate}{\change{Let us first consider a discrete time model to give a motivation on what is the proper formulation of the continuous time model. Given $M$, the time step $\Delta t$ is $\Delta t=\frac{T}{M}$, and denote by $t_k:=k\Delta t ,k=0,1,...,M$ the discrete time points. Consider a purely exploring policy that the agent chooses to stop at time $t_k$ with probability $p^M$ if he haven't stopped earlier. Denote by $\tau^M$ the random time that the agent stops. Then, it is easy to see that 
$$
P(\tau^M>t)=(1-p^M)^{\lfloor \frac{t}{\Delta t}\rfloor}.
$$  
The proper choice of $p^M$ should make the random time distributed over the whole interval $[0,T]$ when $\Delta t \rightarrow 0$. More precisely, we require that, for any $t \in (0,T)$, 
$$
\liminf_{M \rightarrow \infty} P(\tau^M>t) >0 \text{ and } \liminf_{M \rightarrow \infty} P(\tau^M\leqslant t) >0
$$ 
This implies that $(1-p^M)^{\lfloor \frac{t}{\Delta t}\rfloor}=O(1)$, which means that $p^M=O(\Delta t)$. Otherwise, the agent will stop  very early or never stop with probability close to $1$, when $M$ is sufficiently large. Thus, we may write $p^M=\pi \Delta t$. We borrow the terminology from Cox process (see \cite{jeanblanc2009mathematical}) and call $\pi$ the intensity.

The other ingredient for randomized problems in RL is the entropy in the reward functional. As our policy is binary distributed, it is $H(p^M):=-p^M\log p^M-(1-p^M)\log (1-p^M)$, which takes the maximum at $p^M=1/2$. Then, when one adds it into the reward functional, this will prevent us to choose $p^M$ with the right scale. To see this, one can compute that 
\begin{equation*}
	\begin{split}
		H(\pi \Delta t)&=-\pi \Delta t \log{\pi \Delta t }-(1-\pi \Delta t) \log {(1-\pi \Delta t)}\\
		&=-\pi\Delta t\log \Delta t +(\pi-\pi \log \pi)\Delta t+ o(\Delta t).
	\end{split}
\end{equation*}
As the leading term is $\Delta t\log \Delta t$, one should choose $\pi$ as large as possible to maximize $H(p^M)$. This will make the agent choose to stop too early, which contradicts to what we want by exploring. For this reason, we need to adjust the regularizer in the reward functional. Instead of entropy on the probability $p^M$, we need a function defined on the intensity $\pi$. From previous analysis, a natural choice is the $\Delta t$ term, i.e. $R(\pi):=\pi-\pi\log\pi$. We will see that, with this regularizer, the problem is mathematical tractable and the optimal policy has a simple, analytical form. Note that $-R(\pi)=\pi\log \pi-\pi $ is called unnormalized negentropy and used in optimization theory, see \cite{zimmert2019connections}.  Adding a constant, we see that $\pi\log\pi-\pi+1=D(\pi,1)$, where $D(X,Y)=\text{tr}(X\log X-X\log Y-X+Y)$ is the quantum relative entropy defined for positive definite matrices, see \cite{kulis2009low} and references therein.
}}
\subsection{Exploratory Formulation}
\Copy{MathRigor1}{
\change{
Now, let us give the continuous time model. We introduce an additional state $p_t$, which represents the conditional probability that the agent do not stop before time $t$. Motivated by previous discussion, we assume that its dynamic follows
$$
 dp_t=-\pi_t p_t dt, p_0=1.
$$
Here $\pi_t$ is called intensity and is our control variate.  One can treat $\pi_t d t$ as the probability that the agent choose to stop between $[t,t+d t]$ conditioning on he hasn't stopped before $t$ and $\pi_tp_td t$ is the unconditioned probability that the agent stop between $[t,t+d t]$. Then, the reward the agent received will be  
$$
\mathbb E\left[  \int_0^T e^{-rt}g(S_t)\pi_t p_t dt+e^{-rT}g(S_T)p_T \right]
$$
with $g(S)=(K-S)^+$. Note that Touzi and Vieille \cite{touzi2002continuous} also study randomized stopping time. Their main idea is to identify stopping times with $\{0,1\}$-valued, nondecreasing processes. Then convexifying the set of these processes leads naturally to considering the set of all adapted, nondecreasing, right-continuous processes  $p$ with $p_{0-}=1$ and $p_T\le 1$. In our formulation, we actually consider a more restricted set of all adpated, nondecreasing processes with absolute continuous paths.

In additionally, the agent will receive reward by the regularizer before he stops. As we have discussed, it is more reasonable to use the function $R(\pi)$ on intensity than the entropy on probability. \Copy{r2c4}{Thus, the final reward functional is 
\begin{equation}\label{reward_func}
	\mathbb E\left[  \int_0^T e^{-rt}g(S_t)\pi_t p_t+\lambda e^{-rt} R(\pi_t)p_t dt+e^{-rT}g(S_T)p_T \right]
\end{equation}
with $R(\pi)=\pi-\pi\log \pi$ and $\lambda>0$ being a constant.}  $\lambda$ is an exogenous exploration weight parameter capturing the trade-off between exploitation  and exploration. It is usually called temperature constant.  

To sum up, our optimal control problem is formulated as follows.  Let $(\Omega,\mathcal F,\mathbb P)$ be a probability space, on which a standard Brownian motion $W$ is defined. $\mathbb F$ is the filtration generated by $W$ that satisfies the usual condition. The admissible set $\mathcal A$ is defined as the totality of all $\mathbb F$-adapted non-negative processes. For any $t$, given $\pi \in \mathcal A$, $p \in (0,1)$ and $S \in \mathbb R^+$, the controlled dynamic is 
\begin{equation}
	\left\{
	\begin{split}
		dp_s&=-\pi_s p_s ds, p_t=p\\
		dS_s&=S_s\left[r ds+\sigma dW_s \right],S_t=S.
	\end{split}
	\right.
\end{equation} 
As usual, taking conditional expection conditioning on initial state, the reward functional is 
\begin{equation}
	\small
	\begin{split}
	&J(S,p,t;\pi)\\
	=&\mathbb E\left[  \int_t^T e^{-r(s-t)}g(S_s)\pi_s p_s+\lambda e^{-r(s-t)} R(\pi_t)p_s dt+e^{-r(T-t)}g(S_T)p_T\bigg|p_t=p,S_t=S \right],
	\end{split}
\end{equation}
Furthermore, the value function is defined as
$$
\tilde V(S,p,t)=\sup_{\pi \in \mathcal A}J(S,p,t;\pi),
$$
}}

\subsection{Probabilistic Formulation}
In this subsection, we introduce an equivalent formulation for our model. Let $\Theta$ be a random variable, which is exponential distributed and independent from Brownian motion. Given $\{ \pi_t\}_t$, define a random time $\tau$ as 
$$
\tau:= \inf \left\{ t\in [0,T]: \int_0^t \pi_s ds \ge \Theta\right\}  {\bigwedge} T,
$$
where we adopt the convention that the infimum of an empty set is infinity. This construction is referred as the Cox process model, which has been extensively studied in credit risk modelling. The readers are referred to \cite{jeanblanc2009mathematical} for details. One key feature  is that, given the $\sigma$-algebra $\mathcal F_t$, the conditional distribution function of $\tau$ is 
\begin{equation}\label{cond_dist}
P(\tau>s|{\mathcal F}_t)=\exp(-\int_0^s \pi_udu)=p_s,
\end{equation}
for $s \le t$. Moreover, from \cite[Lemma 7.3.4.3]{jeanblanc2009mathematical}, if $h$ is an $\mathbb F$-predictable bounded process, then 
$$
\mathbb E\left[h_{\tau\wedge t}\right]=\mathbb E\left[ \int_0^T h_u \pi_u p_u du+h_Tp_T\right].
$$
Hence, we see that the reward functional \eqref{reward_func} is equivalent to 
$$
\mathbb E\left[ e^{-r\tau}g(S_{\tau }) +\lambda \int_0^\tau e^{-rt}R(\pi_t)dt\right].
$$
From this formulation, we can treat the random time $\tau$ as the execution time.  Note that $\tau$ is not a $\mathbb F$-stopping time. Hence, given $\mathcal F_t$, one can not know whether or not the execution has happened.  In fact, conditioned on that the execution does not take place  before  time $t$ and given  $\mathcal F_t$, the probability that $\tau \in [t,t+\Delta t]$ is $\pi_t \Delta t$. In this point of view,  our model resembles to the exploratory framework of \cite{wang2020reinforcement} that the investors choose to act randomly. 
\section{HJB equation for the related problem}\label{sect_analy}
It is straightforward to derive that the optimal value function $\tilde V(S,p,t)$ satisfies the following HJB equation
$$
\partial_t \tilde V +\mathcal L_S \tilde V+\sup_{\pi \in[0,\infty)} -\pi p \partial_p \tilde V+g(S)\pi p+\lambda R(\pi)p =0,\tilde V(S,p,T)=g(S)p.
$$ 
with $\mathcal L_s =\frac{1}{2}\sigma^2 S^2 \partial_{SS}+rS\partial_{S}-r$. We make the ansatz that $\tilde V(S,p,t)=pV(S,t)$. Then, we have 
\begin{equation}\label{hjb_1}
	\partial_t V+\mathcal L_S  V+\sup_{\pi \in[0,\infty)} (g-V)\pi +\lambda R(\pi) =0,V(S,T)=g(S).
	 \end{equation}
To emphasize its dependence on $\lambda$, we denote by  $V^{\lambda}$ as the solution of \eqref{hjb_1} and still call it the optimal value function. From the optimality condition, the optimal intensity should be  $\bar \pi^\lambda:=\exp(-\frac{V^{\lambda}-g}{\lambda})$. Then, we can rewrite \eqref{hjb_1} as 
\begin{equation}\label{hjb_2}
	\partial_t V^{\lambda}+\mathcal L_S  V^{\lambda}+\lambda \exp(-\frac{V^{\lambda}-g}{\lambda}) =0,V^{\lambda}(S,T)=g(S). 	
\end{equation}
Finally, let $S=e^x$ and $u^{\lambda}(x,t)=V^{\lambda}(S,t)$. It is easy to see that $u^{\lambda}$ solves
\begin{equation}\label{lambda_HJB_2}
	\partial_t u^{\lambda}+\mathcal L_x  u^{\lambda}+\lambda \exp(-\frac{u^{\lambda}-h}{\lambda}) =0,u^{\lambda}(x,T)=h(x),
\end{equation}
with $\mathcal L_x=\frac{1}{2}\sigma^2 \partial_{xx} +(r-\frac{1}{2}\sigma^2)\partial_x-r$ and $h(x)=g(e^x)$. For this equation, we shall have the following theorem for the solvability.
\Copy{MathRigor2thm}{
\begin{thm}\label{thm_solvability}For any $\lambda>0$, there exists a unique classical solution  $u^{\lambda} \in C^{2+\alpha,1+\alpha/2}_{loc}(\mathbb R \times [0,T)) \cap C(\mathbb R \times [0,T])$ for \eqref{lambda_HJB_2} with any $0<\alpha<1$. It holds that
	\begin{equation}\label{esti_bound}
		0\le u^{\lambda}(t,x) \le  K+\lambda(T-t).
	\end{equation} 
\change{ Moreover, we have that the optimal value function $\tilde V(S,p,t)$ equals $u^\lambda(log S,t)p$ and the optimal control $\bar \pi^\lambda$ is given by 
	$$
	\bar \pi^\lambda_t=\exp(-\frac{u^\lambda(\log S_t,t)-g(S_t)}{\lambda}).
	$$

}
\end{thm}}
\begin{rmk}\label{rmk_1}In this paper, we adopt the terminology from RL literature and call the value function of an optimal control problem as optimal value function. Given a feedback strategy  $\pi$,  the expected reward functional under this strategy is called  the value function and denoted as $V^{\pi}(p,S,t)$.  We can make a similar ansatz as before that $V^{\pi}(p,S,t)=pu^{\pi}(\log S,t)$ with $u^{\pi}(x,t)$ being a function defined on $\mathbb R \times [0,T]$. With a little bit of abuse of terminology, we also call $u^{\pi}$ the value function.   	
\end{rmk}
\subsection{Policy Iteration}
In reinforcement learning, one method to learn the optimal strategy is policy iteration in which the optimal strategy is approximated by iteratively updating. More precisely, given a feedback  strategy $\pi^n(x,t)$, the corresponding value function $u^n(x,t)$ (see Remark \ref{rmk_1}) satisfies 
\begin{equation} \label{eq_VI}
\partial_t u^n+\mathcal L_x u^n+H(x,\pi^n(x,t),u^n)=0,u^n(x,T)=h(x),
\end{equation}
with the Hamiltonian $H$ being defined as
$$
H(x,\pi,u)=(h(x)-u)\pi+\lambda R(\pi).
$$
Having a value function $u^n$, one can construct a feedback strategy $\pi^{n+1}$ as 
\begin{equation}\label{eq_Policy_Update}
	\pi^{n+1}(x,t)=\mathop{\text{argmax}}_{\pi \in(0,\infty)}H(x,\pi,u^n(x,t))=\exp(-\frac{u^n(x,t)-h(x)}{\lambda})
\end{equation}
We continue this iteration and obtain a sequence of pairs of strategy and value function. The following theorem states that the performance of the policies is improved during the iteration.
\begin{thm}\label{thm_monotone}
	Give any initial guess $u^0 \in C_b(\mathbb R)$ for the optimal value function. $\{u^n,\pi^n\}_{n=1,2,...}$ are defined iteratively according to \eqref{eq_Policy_Update} and \eqref{eq_VI}. Then, we have that 
	$u^{\lambda}\ge u^{n+1} \ge u^{n}$, for $n=1,2,...$.
\end{thm} 
One should expect that the value functions $\{u^n\}$ converge to the optimal value function $u^{\lambda}$. This is indeed the case. Furthermore, we can give an estimate for the convergence rate of value functions.
\begin{thm} \label{thm_convergence_rate} Define by $M:=\|(h-u^1)^+\|_{\infty}$. There exists a constant $C$ independent of $\lambda$ such that 
	$\|u^{n+1}-u^\lambda\|_{\infty} \le C^n\frac{T^n}{n!}\exp(n(\frac{M}{\lambda}+rT))(M+K+\lambda T)$.
\end{thm}
\Copy{MathProof1}{
\change{Due to the Stirling formula, it holds that $n!\backsim  \exp(n\log n) n \sqrt{\frac{2\pi}{n}} e^{-n}$.  Hence, we see that the right hand side of the estimation in Theorem \ref{thm_convergence_rate} will go to zero as $n$ goes to $\infty$, no matter what the values of $C,T,M,r$ are. This gives the convergence of policy iteration.}}
\begin{rmk}
\begin{enumerate}[i)]
	\item In \cite{kerimkulov2020exponential}, the authors proved an exponential convergence rate $Cq^n$ for policy iteration. Our rate here is faster than exponential with respect to $n$. We think that the main reason for this is that the optimal $
	\pi$  smoothly depends on other variables of the Hamiltonian $H$, see \eqref{eq_Policy_Update}. For general optimal control problem, such a property does not hold.
	\item  Note that $u^1$ is defined by \eqref{eq_VI}. Thus, $M$ should also depend on $\lambda$. In fact, one can get an estimation of $M$, which is decreasing with respect to $\lambda$.
 But, for the special initial value function $u^0=g$, $M$ is independent of $\lambda$. 
\end{enumerate}
\end{rmk}
\subsection{Comparison between classical problem and randomized problem}
For the classical American put option pricing problem, it is well-known that the option price is given by $u^*(\log S,t)$ with the function $u^*$ being the solution of the following variational inequality
\begin{equation}\label{VI}
\min\{-\partial_t u^*-\mathcal L_x u^*,  u^*(x,t)-h(x) \}=0,u(x,T)=h(x).
\end{equation}
and the optimal stopping time $\tau^*$ is defined as 
\begin{equation*}
	\tau^*:=\inf\{ t| u^*(\log S_t,t) \le g(S_t)\} \wedge T.
\end{equation*}
It is easy to see that \eqref{VI} is formally the limit of \eqref{hjb_2} as $\lambda$ going to $0$. In other words, \eqref{hjb_2} can be seen as a penalized equation for \eqref{VI}. Thus, it is reasonable to believe that $u^{\lambda}$ should converge to $u^*$ as $\lambda$ going to $0$. To prove this,  the following lemma is needed. \change{
\begin{lemma}\label{lem_lower_bound}
	We have, for some constant $C$ independent of $\lambda$, 
	$$u^{\lambda}(x,t)\ge h(x)-\lambda \log\frac{CrK}{\lambda}.$$
\end{lemma}}
Then, we have the following theorem.\Copy{r2c41}{
\begin{thm}\label{thm_difference}
	We have
	$$
	u^\lambda(x,t)-\lambda \log(1+(T-t))\le u^*(x,t) \le u^{\lambda}(x,t)+ \lambda   \log \frac{CrK}{\lambda},
	$$
	where $C$ is the constant in Lemma \ref{lem_lower_bound}.
\end{thm}}

\Copy{MathRigor3}{
Having the value function $u^{\lambda}$ or $V^{\lambda}$ equivalently, how can we construct a strategy for the optimal stopping problem? Maybe the most straight forward way is to execute at the stopping time 
\begin{equation}\label{stopping_rule}
	\tau^{\lambda}:=\inf\{ t| V^{\lambda}(S_t,t) \le g(S_t)\} \wedge T.
\end{equation}
\change{Clearly, $\tau^{\lambda}$ is sub-optimal for the American option problem.  Thus, it is obvious that 
	$\mathbb E\left[e^{-r\tau^\lambda}g(S_{\tau^{\lambda}})  \right] \le \mathbb E\left[e^{-r\tau^*}g(S_{\tau^*}) \right]$. However, note that the optimal stopping time $\tau^*$ is defined in a similar way as $\tau^\lambda$, but replacing $V^\lambda$ in \eqref{stopping_rule} by $V^*$. We have shown that $V^\lambda$ converges to $V^*$ in Theorem \ref{thm_difference} as $\lambda$ goes to zero. Then, one should expect that   $\mathbb E\left[e^{-r\tau^\lambda}g(S_{\tau^{\lambda}})  \right]$ is close to  $\mathbb E\left[e^{-r\tau^*}g(S_{\tau^*}) \right]$. The following theorem confirms this assertion and further proves that the difference is at most on the scale of $\lambda \log \lambda$.
} 
\Copy{r2c42}{ \begin{thm}\label{thm_wealth_loss}
 	$\mathbb E\left[e^{-r\tau^\lambda}g(S_{\tau^{\lambda}})  \right] \ge \mathbb E\left[e^{-r\tau^*} g(S_{\tau^*}) \right]-O(\lambda \log\frac1\lambda)$.
 \end{thm}}}
\begin{rmk}
	\begin{enumerate}[i)]
		\item Recently, Tang et al \cite{tang2021exploratory} studied the exploratory HJB equation arising from the entropy-regularized exploratory control problem. For some special case, they also proved an convergence rate of $O(-\lambda \log \lambda)$ for the value functions as $\lambda \rightarrow 0^+$.
		\item We show that the convergence rate is controlled by  $\frac{(CT)^n}{n!} \exp(\frac{nM}{\lambda})(M+K+\lambda T)$ which is decreasing in $\lambda$ when $\lambda$ is small. On the other hand, the bias is dominated by $O(-\lambda \log(\lambda))$. This implies a balance between convergence speed and bias in the choice of $\lambda$. 
	\end{enumerate}
	
\end{rmk}
\subsection{Numerical Solution for \eqref{lambda_HJB_2}}In this subsection, we numerically solve  \eqref{lambda_HJB_2} for illustration. To deal with the nonlinear term, we adopt the method in \cite{forsyth2002quadratic} which combines  finite-difference method with generalized Newton iteration. The grid of mesh points is taken to be $(t_i,x_j)=(i\Delta t,j\Delta x)$ with $i=0,1,...,M$ and $j=0,\pm 1,...,\pm N$, where $M=\frac{T}{\Delta t}$. Let $U_{i,j}$ be the discrete solution at mesh point $(t_i,x_j)$ and $A$ be the discretization of the second order operator $-\mathcal L_x$. The iteration is as follows
\begin{breakablealgorithm} \caption{Finite-difference Iterative Algorithm}
	\begin{algorithmic} 	
		\Require{$\Delta t,\Delta x$, $\lambda$ and tolerance $tol$}
		\State{Compute $A$}
		\State{Initialize $U_{M,j}=(K-e^{x_j})^+$}
		\For{ $t=M-1,M-2,...,1,0$}
		\State{Set $U^0_{\cdot}=U_{t-1,\cdot}$}
		\For{$k=1,2,3,...$}
		\State{Solve the linear equation
			\begin{equation*}
				\frac{1}{\Delta t}(U^k_j-U_{t+1,j})+[AU^k_{\cdot}]_j=\lambda \exp(-\frac{U^{k-1}_j-U_{0,j}}{\lambda})+\exp(-\frac{U^{k-1}_j-U_{0,j}}{\lambda})(U^{k-1}_j-U^k_j)
			\end{equation*}
		for $j=0,\pm 1,...,\pm(N-1)$, and
		$$
		U^k_{\pm N}=(K-e^{x_{\pm N}})^+.
		$$ 
		} 
		\State{If $\frac{\|U^{k}_{\cdot}-U^{k-1}_{\cdot}\|_{\infty}}{\max\{1,\|U^{k-1}_{\cdot}\|_{\infty}\}}<tol$, Quit}
		\EndFor
		\State{Set $U_{t,\cdot}=U^k_{\cdot}$.}
		\EndFor
	\end{algorithmic}
\end{breakablealgorithm}

\begin{figure}[h]
	\centering
	\includegraphics[width=0.60\textwidth]{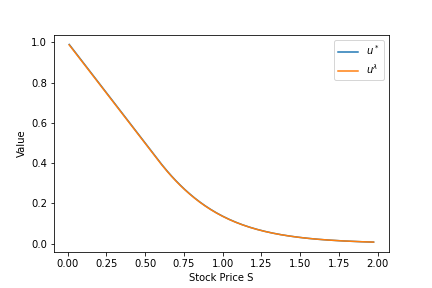}
	\caption{Value functions $u^\lambda$ and $u^*$ at $t=0$}
	\label{VF}
\end{figure}

\begin{figure}[h]
	\centering
	\begin{subfigure}{0.45\textwidth}
		\centering
		\includegraphics[width=1\textwidth]{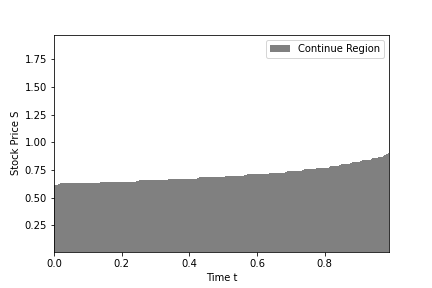}
		\subcaption{ $\{u^\lambda > h\}$}
	\end{subfigure}
	\begin{subfigure}{0.45\textwidth}
		\centering
		\includegraphics[width=1\textwidth]{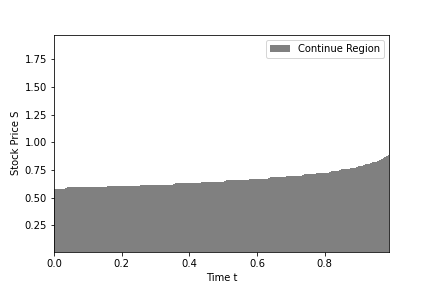}
		\subcaption{ $\{u^* > h\}$}
	\end{subfigure}
	\caption{Continue regions for two functions }
	\label{fig:exection region}
\end{figure}

\begin{figure}[h]
	\centering
	\includegraphics[width=0.60\textwidth]{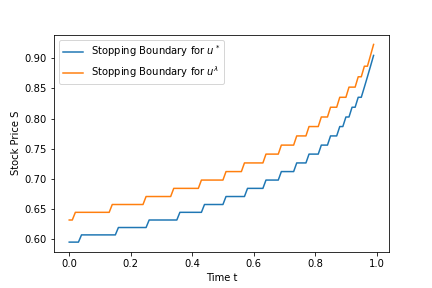}
	\caption{Comparison of the stopping boundaries}
	\label{BD}
\end{figure}

For numerical experiment, we set $r=5\%,\sigma=0.4,k=1,T=1$ and $\lambda=0.005$. In Fig \ref{VF}, we plot the value function $u^\lambda$ at time $t=0$ and compare it with $u^*$, which is computed by using penalty (see \cite{forsyth2002quadratic}). The difference between two functions is very small. In fact, numerical result shows that the maximum error  is about $0.0023$ which is of the same scale as $-\lambda \log \lambda=0.0038$. 
In Figure \ref{fig:exection region}, we plot continue regions for  $u^\lambda$  and $u^*$. It is easy to see that the regions are of similar shape and both separated by stopping boundary.   
For further comparison, we plot these two boundaries in Figure \ref{BD}. One can observe that the stopping boundary of $u^*$ is  lower than that of $u^\lambda$.  This is probably due to the fact that the investor is compensated by the regularization term and thus tends to be continue when the price is relatively higher.
\section{Reinforcement Learning Algorithm}\label{sect_rl} In this section, we propose an reinforcement learning algorithm for the optimal stopping problem. Note that our formulation transforms the optimal stopping problem into a standard optimal control problem. Thus, we could use well-developed reinforcement learning algorithms to learn the execution strategy from the data. In this paper, we adopt an algorithm based on policy iteration. It consists of two  iterating steps: policy  evaluation and policy update. 

\subsection{Policy Iteration}
In numerical implementation, we discretize the time interval $[0,T]$ into $L$ points: $0=t_0<t_1<...<t_L=T$ with equal distance $\Delta t$. In policy iteration step, our task is to calculate the value function $V^\pi$ for a given feedback policy $\pi=\pi(S,t)$. It is clear that $V^\pi(S,T)=g(S)$. At other point $t<T$, it is characterized by the PDE:
$$
\partial_t V^\pi+ \mathcal L_S V^\pi +(g-V^\pi) \pi+\lambda R(\pi)=0.
$$
One can not direct compute this equation as we assume that the model information is not available. Instead, we shall use a parameterized function $V_\theta(\cdot)$ to approximate the value function $V^\pi(\cdot,t)$ with a data set. Denote by $\theta^l$ the corresponding parameters at time $t_l$. For reinforcement learning algorithm, the available data are  samples of the stock prices $\{ S^m_{t_l}\}_{l=0,1,...,L,m=1,2,...,M}$. \Copy{MathRigor4}{ Given a policy $\pi=\pi(S,t)$, one can calculate the state $p$ at time point $t_l$ iteratively by the following rule:
\begin{equation}\label{discrete_p}
	p^m_{t_{l+1}}=\left(1-\pi(t_l,S^m_{t_l})\Delta t\right) p^m_{t_l},p^m_0=1.
\end{equation}
From It\^o-formula, we have that 
\begin{equation}\label{eq_ito}
	\begin{split}
	&\mathbb E\left[ \int_{t_l}^{t_{l+1}} e^{-ru}g(S_u)\pi_up_u+\lambda e^{-ru} R(\pi_u)p_u du+e^{-rt_{l+1}}V^\pi(S_{t_{l+1}},t_{l+1})p_{t_{l+1}}\big|\mathcal F_{t_l} \right]\\
	=&e^{-rt_l}V^\pi(S_{t_l},t_l)p_{t_l}.
	\end{split}
\end{equation}
Recall that conditional expectation can be treated as an $L^2$-projection in  certain Hilbert space.
This motivates us to choose $\theta^l$ that minimizes the following quantity
\begin{equation}\label{eq_td_error_prev}
	\sum^M_{m=1}\left( g(S^m_{t_l})\pi^m_{t_l}p^m_{t_l}\Delta t+\lambda R(\pi^m_{t_l})p^m_{t_l}\Delta t +e^{-r\Delta}V_{\theta^{l+1}}(S^m_{t_{l+1}})p^m_{t_{l+1}}   -V_{\theta^{l}}(S^m_{t_l})p^m_{t_l}\right)^2.
\end{equation}
\change{From \eqref{discrete_p}, we subtract  $p_{t_l}$ in \eqref{eq_td_error_prev}
	and propose minimization of a similar term}
\begin{equation}\label{eq_td_error}
\delta_l:=\sum_{m=1}^M\left( g(S^m_{t_l})\pi^m_{t_l}\Delta t+\lambda R(\pi^m_{t_l})\Delta t +e^{-r\Delta}V_{\theta^{l+1}}(S^m_{t_{l+1}})(1-\pi^m_{t_l}\Delta t)  -V_{\theta^{l}}(S^m_{t_l})\right)^2,	
\end{equation}
which is simpler than previous without the need to compute $p_{t_l}$.} 
 This term is also called TD-error in the context of RL, see \cite{sutton2011reinforcement}. During each iteration, one can update $\theta^l$ by taking one gradient descent step.  
 
  In the policy update step, having calculate value function $V^\pi$, our task is to give an improvement of $\pi$. In general, this can be obtained by the optimality condition in the HJB equation, which involves another optimization problem. However, in our problem, we see that the minimization problem can be explicitly solved as in \eqref{eq_Policy_Update}. This greatly simplifies the policy update procedure. 

To summarize, we propose the following algorithm
\begin{breakablealgorithm}	\caption{Policy Iteration for Optimal Stopping Problem}
	\begin{algorithmic}
	\Require Number of Iterations $N$, Batch Size $M$, Learning Rate $lr$
	\For {$n=1,2,...,N$} 
	\State Generate a batch of stock price $\{ S^m_{t_l}\}_{l=0,1,...,L,m=1,2,...,M}$ with batch size $M$
	\For{$t=L-1,L-2,...,0$} 
	\State Calculate $\pi^m_{t_l}=\exp(-\frac{V_{\theta^l}(S^m_{t_l})-g(S^m_{t_l})}{\lambda})$
	\State Calculate the TD-error $\delta_h$ according to \eqref{eq_td_error}
    \State Update $\theta_l$ by taking one gradient descent step on $\delta_l$ with learning rate $lr$
	\EndFor
	\EndFor
	\end{algorithmic}
\end{breakablealgorithm}
\quad
\\

\Copy{r2c5}{\change{ The algorithm we design is offline, i.e. it utilizes previously collected data, without additional online data collection. The value function is updated backward in time. Compared with the well-known Longstaff-Schwarz method \cite{longstaff2001valuing} for American option pricing, our method is different in several aspects. First, in Longstaff-Schwarz method, the continuation value is approximated by a linear combination of basis functions, while we use neural networks for approximation. Neural networks can handle large-scale problems with high-dimensional input spaces. This enables our method to deal with high-dimensional cases. Second,  Longstaff-Schwarz method  determines the optimal exercise policy by comparing the option's payoff at each exercise opportunity with the expected future value of holding the option. While in our method, the policy is randomized and the stopping probability is determined by the difference between the payoff and the value function. Although the theoretical optimal strategy is deterministic, stochastic strategy is preferable in learning due to its potential benefits including robustness, improved convergence and imitation learning, see \cite{haarnoja2017reinforcement}. Quantify the difference between Longstaff-Schwarz method and ours in convergence rate and sample efficiency is an interesting problem and is one of our future research topic. Finally, since we reduce the optimal stopping problem to a standard optimal control problem, one can use various RL methods to obtain optimal execution policy. In particular, the algorithm can be online, i.e.  the agent learns and improves its behavior while directly interacting with an environment in real-time. 
}}
\subsection{Numerical Examples}In this subsection, we shall present some numerical examples for our RL algorithm. Although our theoretical analysis is based on American put option model, the algorithm is applicable for other models.
\subsubsection{American Put Option}
\begin{figure}[h]
	\centering
	\includegraphics[width=0.60\textwidth]{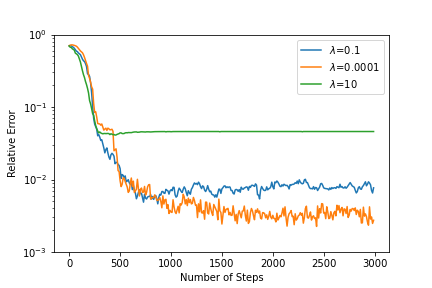}
	\caption{Learning Curves for different $\lambda$}
	\label{lc}
\end{figure}

We consider a standard American put option with $s_0=K=40,r=6\%,\sigma=40\%$. \change{The price has been computed by Becker et al \cite{becker2021solving} as $5.311$ and also by M. Smirnov \cite{smirnov} as $5.318$ using binomial tree method.} We use previous RL algorithm to compute the price. Choose $L=50$. Instead of approximating the value function by neural networks directly, we find that it is more efficient to approximate it by the payoff $g$ plus a neural network. The neural network that we employ consists of $2$ hidden layers with ReLU activation function. In additional to stock price, we also take the payoff of the option as  input of the neural network. The width of hidden layers is $21$. The training procedure consists of $N=3000$ steps with a batch size of $M=2^{10}$. We also generate a test data set with size $2^{18}$.  
\Copy{r2c7}{
To compare the effect for different $\lambda$, we train three neural networks with $\lambda=10$,$10^{-1}$,$10^{-4}$ at the same time. The final performance is $5.052,5.273,5.298$ respectively. \change{Every $10$ steps, we test the learned strategy on the test data set to get an estimation $\mathcal P$ of the expected reward. We compare $\mathcal P$ with the price obtained in \cite{becker2021solving}. The relative error is computed as $\frac{|\mathcal P-5.311|}{5.311}$. }} The learning curves are plotted in Figure \ref{lc}.  It clearly shows the trade-off between  convergence speed and bias. With a larger $\lambda$, the learning curve tends to decrease and converge earlier. But, this also leads to a larger bias for the final performance. 
\Copy{r2c6}{
\change{To see whether the error will continue to decrease as $\lambda$ becomes smaller, we further conduct the experiment with $\lambda$ ranging from $10$ to $10^{-6}$. The result is represented in Table \ref{t4}. One can see that, in general, the performance will be improved for smaller $\lambda$. But, it will be hard to further reduce the error when $\lambda$ is less than $10^{-4}$. We think the major reason is that we use a constant size of batch in the training. When $\lambda$ is small, the sampling error is dominant compared with the bias introduced by $\lambda$.  }
\begin{table}[h]\centering
	\caption{Final performance for different $\lambda$}
	\label{t4}
	\begin{tabular}{|c|c|c|c|c|c|c|c|c|}
		\hline
		$\lambda$ &10&1&$10^{-1}$&$10^{-2}$&$10^{-3}$&$10^{-4}$& $10^{-5}$&$10^{-6}$\\
		\hline
		Performance&5.052&5.131&5.273&5.294&5.291&5.298&5.294&5.295\\
		\hline
	\end{tabular}
\end{table}
}
\subsubsection{Bermudan Max-Call Option}  Bermudan max-call are one of the most studied examples in the numerical literature on optimal stopping problem, see, for example, \cite{longstaff2001valuing} and references therein. The payment is decided by the maximum of $d$ underlying assets. More precisely, the risk neutral dynamics of the assets are characterized by a multi-dimensional Black-Scholes model:
$$
S^i_t=S_0^i\exp(\big(r-\delta_i-\frac{\sigma^2_i}{2}\big)t+\sigma^i W^i_t), i=1,2,...,d,
$$
for initial prices $S_0^i \in(0,\infty)$, risk free rate $r>0$, dividend rates $\delta_i>0$, volatilities $\sigma_i>0$ and a $d$-dimensional Brownian motion with constant correlations $\rho_{ij}$ between components $W^i$ and $W^j$. A  Bermudan max-call option on $S^1,S^2,...,S^d$ has payoff $(\max_i S^i_t-K)^+$. In the numerical examples, the option can be executed by at any time of a time grid $0=t_0<t_1<...<t_N=T$. Then, the price of the option can be give as 
$$
\sup_{\tau}\mathbb E\left[e^{-r\tau}  (\max_i S^i_\tau-K)^+ \right].
$$

We employ the same architecture of the neural network as in American put option case. The width of the hidden layer is taken to be $d+20$. The training procedure consists of $N=5000$ steps with a batch size $M=2^{10}$. For the case $d=500$, we need $N=8000$ steps. The test data set is of the size $2^{18}.$ We set the temperature constant $\lambda=10^{-2}$. For the coefficients of the model, we consider two cases and compare the results with that in \cite{becker2019deep}.

\textbf{Symmetric Case}\\
We first consider a symmetric case, where $S^i_0=s_0,\delta_i=\delta,\sigma_i=\sigma$ for $i=1,2,...,d$ and $\rho_{ij}=\rho$. The results are represented in Table \ref{t1}.

\textbf{Asymmetric Case}\\
We next consider an asymmetric case. All the coefficients except $\sigma_i$ are the same as in symmetric case.  For $d<5$, we set $\sigma_i=0.08+0.32\times(i-1)/(d-1),i=1,2...,d$. For $d>5$, $\sigma_i=0.1+i/(2d),i=1,2,...,d$. The results are represented in Table \ref{t2}.
\Copy{r2c81}{
 \begin{table}[!htbp]\centering
  	\small
  	\caption{Summary result for Bermudan max-Call option with $r=5\%,\delta=10\%,\sigma=20\%,\rho=0,K=100,N=9$.}
	\begin{subtable}{0.75\textwidth}
		\begin{tabular}{|c|c|c|c|c|}
			\hline
			$d$&$s_0$&Our Method&Price in \cite{becker2019deep}&Runtime in sec. \\
			\hline
			2&90&8.067&8.074&196\\
			2&100&13.871&13.903&196\\
			2&110&21.352&21.349&192\\
			5&90&16.646&16.644&211\\
			5&100&26.161&26.159&213\\
			5&110&36.720&36.772&210\\
			10&90&26.240&26.240&220\\
			10&100&38.255&38.337&218\\
			10&110&50.808&50.886&224\\
			50&90&53.983&54.057&473\\
			50&100&69.528&69.736&480\\
			50&110&85.185&85.463&485\\
			100&90&66.299&66.556&829\\
			100&100&83.451&83.584&839\\
			100&110&100.630&100.663&843\\
			200&90&78.936&79.174&2158\\
			200&100&97.459&97.612&2048\\
			200&110&115.506&116.088&2051\\
			500&90&96.322&96.147&6875\\
			500&100&116.327&116.425&6934\\
			500&110&136.849&136.765&7084\\
			\hline
		\end{tabular}
		\subcaption{Symmetric Case}
		\label{t1}
	\end{subtable}

\begin{subtable}{0.75\textwidth}
		\begin{tabular}{|c|c|c|c|c|}
			\hline
			$d$&$s_0$&Our Method&Price in \cite{becker2019deep}&Runtime in sec.\\
			\hline
			2&90&14.025&14.339&193\\
			2&100&19.792&19.808&198\\
			2&110&27.529&27.158&197\\
			5&90&27.579&27.662&210\\
			5&100&37.974&37.985&209\\
			5&110&49.595&49.499&205\\
			10&90&85.997&85.987&219\\
			10&100&104.438&104.741&220\\
			10&110&122.962&123.745&221\\
			50&90&196.286&196.321&490\\
			50&100&277.513&277.886&481\\
			50&110&259.109&259.321&479\\
			100&90&263.522&263.579&844\\
			100&100&302.364&302.765&831\\
			100&110&341.271&341.575&835\\
			200&90&344.165&345.146&2276\\
			200&100&392.764&392.958&2092\\
			200&110&439.654&440.815&2201\\
			500&90&475.867&477.102&7156\\
			500&100&538.678&539.577&7034\\
			500&110&601.567&602.252&6923\\
			\hline
		\end{tabular}
		\subcaption{Asymmetric Case}
		\label{t2}
	\end{subtable}
\end{table}}
\subsubsection{Optimal Stopping for  fractional Brownian Motions}
In this part, we consider a optimal stopping for fractional Brownian motions. A fractional Brownian motion with Hurst parameter $H\in(0,1]$ is a continuous centered Gaussian process $(W^H_t)_{t\ge0}$ with covariance
$$
\mathbb E\left[ W^H_tW^H_s\right]=\frac12\left(|t|^H+|s|^H-|s-t|^H  \right).
$$
For $H=1/2$, $W^H$ is a standard Brownian motion. But, when $H \neq \frac12$, it is neither a martingale nor a Markov process.  We want ot approximate the supreme 
$$
\sup_{0\le\tau\le1}\mathbb E\left[ W^H_\tau \right],
$$
over all $W^H$-stopping time $\tau$. 

Denote by $t_n=n/100,n=0,1,...,100$. As the model is non-Markovian in general, the value function $V$ should be path-dependent, i.e. $V(\cdot,t_n)$ not only depends on current state $W_{t_n}$, but also the whole historical path.  For this reason, we approximate the value function $V(\cdot,t_n)$ by a neural network with $n$-dimensional input. It also consists of two hidden layer with width $n+20$. The training procedure consists of $N=6000$ steps with a batch size $M=2^{11}$. The test data set is of the size $2^{15}$. 
The temperature constant $\lambda$ is chosen to be $10^{-5}$. We consider different values of $H$ and these results are shown in Table \ref{t3}.
\Copy{r2c82}{
\begin{table}[h]
	\centering
	\caption{Summary result for fraction Brownian motion}
	\label{t3}
\begin{tabular}{|c|c|c|c|}
	\hline
	H&Our Method&Price in \cite{becker2019deep}&Runtime in sec.\\
	\hline
	0.01&1.513&1.519&4636\\
	0.05&1.283&1.293&4669\\
	0.10&1.043&1.049&4558\\
	0.15&0.825&0.839&4587\\
	0.20&0.651&0.658&4700\\
	0.25&0.494&0.503&4461\\
	0.30&0.364&0.370&4540\\
	0.35&0.243&0.255&4567\\
	0.40&0.144&0.156&4678\\
	0.45&0.061&0.071&4568\\
	0.50&0.000&0.002&4632\\
	0.55&0.052&0.061&4661\\
	0.60&0.108&0.117&4553\\
	0.65&0.155&0.164&4626\\
	0.70&0.192&0.207&4541\\
	0.75&0.233&0.244&4714\\
	0.80&0.261&0.277&4650\\
	0.85&0.299&0.308&4705\\
	0.90&0.331&0.337&4645\\
	0.95&0.354&0.366&4667\\
	1.00&0.395&0.395&4578\\
	\hline
	\end{tabular}
\end{table}}
\section{Conclusion} In this paper, we study the randomized optimal stopping time problem in the so-called exploratory framework. By introducing an additional state, we reduce it to a standard optimal control problem. Using PDE techniques, we are able to give some quantitative results, especially on the trade-off between convergence speed and bias. Although we mainly focus on American put option with one risky asset, these theoretical results can be extended to more general Markovian models as long as some regularity conditions hold.  For non-Markovian cases, we believe that there will be similar results. But one may need the theory of backward stochastic differential equations to give a rigorous proof. We also adjust well-developed RL algorithm to learn the optimal stopping strategy. Computational results show that our algorithm achieves comparable performance as that in \cite{becker2019deep,reppen2022neural}.  Noting that the algorithm here is offline, one can also apply an online algorithm to study the strategy.

Our future research is to extend this framework to other optimal control problems in mathematical finance. For optimal switching problem, the extension is straightforward. Instead of introducing an additional scalar-valued state, one needs  a vector-valued state to represent the distribution of the population in different modes, Moreover, the optimal strategy will be a soft-max type. But, it requires more effort to give some rigorous theoretical results. For optimal trading problem with transaction cost and tax, the situation is more complicated as the free boundary is also determined by the derivatives of  the value function.       
\section{Proofs of the main results}\label{sect-proof}
\subsection{Proofs for Theorem \ref{thm_solvability}}\Copy{MathProof3}{
The proof is rather standard in parabolic PDE theory. One has to deal with the unboundedness of the exponential function. Let $\{h_\varepsilon(x)\}$ be a sequence of smooth functions such that approximating $h(x)=(K-e^x)^+$ increasingly and uniformly.  Given $M >0$, consider a smooth cut-off function $\phi_M$ such that $\phi_M(x)=e^x$ for $x \le M$ and $\phi_M(x)=e^{M+1}$ for $x \ge M+1$. Hence, $\phi_M$ is bounded and Lipschitz continuous. Denote by $\Omega_N:=(-N,N) \times [0,T)$. First, we will solve the following PDEs
\begin{equation}
	\left\{
	\begin{split}
		&\partial_t u^{\lambda,N,\varepsilon}+\mathcal L_x  u^{\lambda,N,\varepsilon}+\lambda \phi_M(-\frac{u^{\lambda,N,\varepsilon}-h_\varepsilon}{\lambda}) =0, (x,t) \in \Omega_N,\\
		&u^{\lambda,N,\varepsilon}(-N,t)=h_\varepsilon(-N)+\lambda(T-t),u^{\lambda,N,\varepsilon}(N,t)=h_\varepsilon(N)+\lambda(T-t), t\in[0,T]\\
		&u^{\lambda,N,\varepsilon}(x,T)=h_\varepsilon(x),x \in [-N,N].\\ 
	\end{split}
	\right .
	\label{PDE_cutoff}
\end{equation}
\change{
As $\mathcal L_x$ is a second order operator with constant and non-degenerate coefficients, classical results for parabolic PDEs, including maximum principle, $L^p$ and Schauder estimations hold in our case.
Since the nonlinear term  $(u,x) \mapsto \lambda \phi_M(-\frac{u-h_\varepsilon(x)}{\lambda})$ is bounded and Lipschitz continuous with respect to $(u,x)$, one can apply the existence theorems for nonlinear PDEs (see \cite[Sec 7.4, Theorem 8]{friedman2008partial}) to get that there exists a classical solution $u^{\lambda,N,\varepsilon}$, which belongs to $C^{2+\alpha,1+\alpha/2}(\Omega_N)\cap C(\bar \Omega_N)$ for any $\alpha \in (0,1)$.     

 Noting that $(\partial_t+\mathcal L_x)u^{\lambda,N,\varepsilon}=-\lambda \phi_M(-\frac{u^{\lambda,N,\varepsilon}-h_\varepsilon}{\lambda}) <0$, $u^{\lambda,N,\varepsilon}(\pm N,t) \ge 0$ and $u^{\lambda,N,\varepsilon}(x,T) \ge 0$, classical maximum principle \cite[Theorem 2.4]{lieberman1996second} implies that $u^{\lambda,N,\varepsilon} \ge 0$. }Define $w(x,t):=u^{\lambda,N,\varepsilon}(x,t)-\big(K+\lambda(T-t)\big)$. Then, we see that $w$ satisfies the following equation in $\Omega_N$:
$$
\partial_t w +\mathcal L_x w+\lambda \phi_M(-\frac{u^{\lambda,N,\varepsilon}-h_\varepsilon}{\lambda})=\lambda+rK+r\lambda(T-t),
$$
with $w(x,T)\le 0$. We can rewrite it as 
\begin{equation*}
	\begin{split}
		\partial_t w +\mathcal L_x w+\lambda \phi_M(-\frac{u^{\lambda,N,\varepsilon}-h_\varepsilon}{\lambda})-\lambda \phi_M(-\frac{K+\lambda(T-t)-h_\varepsilon}{\lambda}) \\
		=\lambda - \lambda \phi_M(-\frac{K+\lambda(T-t)-h_\varepsilon}{\lambda})+rK+r\lambda(T-t).
		\end{split}
\end{equation*}
Since $(K+\lambda(T-t))\ge h(x) \ge h_\varepsilon(x)$, we have that the right hand side is non-negative. Using maximum principle again, we have that $w \le 0$, i.e. $u^{\lambda, N,\varepsilon}(x,t) \le K+\lambda(T-t)$. Now, we see that the bound is independent of $M$. Choosing $M$ sufficient large, $u^{\lambda, N,\varepsilon}$ solves \eqref{PDE_cutoff} with $\phi_M(\cdot)$ replaced by $\exp(\cdot)$.

\change{Choose $p>3$. For any $R,\eta>0$ and $N>R+1$, classical interior $L^p$ estimates \cite[Theorem 7.22]{lieberman1996second} yields that,
\begin{equation*}
	\small
	\begin{split}
	&\|u^{\lambda ,N,\varepsilon}\|_{W^{2,1}_p((-R,R)\times(\eta,T))}\\
	 \le& C(R,\eta,p)\left( \|u^{\lambda ,N,\varepsilon}\|_{L^p((-R-1,R+1)\times(0,T))}+\|\lambda\exp(-\frac{u^{\lambda ,N,\varepsilon}-h_\varepsilon}{\lambda})\|_{L^p((-R-1,R+1)\times(0,T))}\right)	
		\end{split}
\end{equation*}
As we have already proved the boundedness of $u^{\lambda ,N,\varepsilon}$, it holds that 
$$
\|u^{\lambda ,N,\varepsilon}\|_{W^{2,1}_p((-R,R)\times(\eta,T))} \le C(\lambda,R,\eta,p).
$$
The parabolic version of embedding theorem (see \cite[Theorem 3.14]{hu2011blow}) implies that $W^{2,1}_p((-R,R)\times(\eta,T))$ can be  embedded into $C^{1+\alpha,\frac{1+\alpha}{2}}([-R,R]\times[\eta,T])$ with $\alpha=1-\frac3p$, i.e.
$$
\|u\|_{C^{1+\alpha,\frac{1+\alpha}{2}}([-R,R]\times[\eta,T]} \le C(R,\eta,T,p) \|u\|_{W^{2,1}_p((-R,R)\times(\eta,T))}.
$$
Hence,  we see that $u^{\lambda ,N,\varepsilon}$ is uniformly equi-continuous on $[-R,R]\times[\eta,T]$.
Then, it  yields that, subtracting  a subsequence if needed, $\{u^{\lambda,N,\varepsilon}\}$ locally uniformly converge and also weakly converge in $W^{2,1}_{p,loc}(\mathbb R\times(0,T))$  to a function $u^{\lambda,\varepsilon}$ as $N$ goes to infinity. One can verify that $u^{\lambda,\varepsilon}$ solves 
$$
\partial_t u^{\lambda,\varepsilon}+\mathcal L_x u^{\lambda,\varepsilon}+\lambda \exp(-\frac{u^{\lambda,\varepsilon}-h_\varepsilon}{\lambda})=0, (x,t)\in \mathbb R \times(0,T).
$$
Let us check that $u^{\lambda,\varepsilon}(x,T)=h_\varepsilon(x)$. For that purpose, consider the following PDE:
\begin{equation*}
	\left\{
	\begin{split}
		&\partial_t v	+\mathcal L_x v=0, (x,t) \in \mathbb R \times (0,T),\\
		&v(x,T)=h_\varepsilon(x).
	\end{split}
	\right.
\end{equation*} 
It is easy to get that $v$ is bounded. Then, the point-wise bound estimate (see \cite[Theorem 2.10]{lieberman1996second}) yields that
\begin{equation*}
	\small
	\begin{split}
		&\sup_{|x|<N } |u^{\lambda,N,\varepsilon}(x,t)-v(x,t)| \\
		\le& e^{r(T-t)}\left( \sup_{|x|<N,s\in[T,t]}|\lambda \exp(-\frac{u^{\lambda,N,\varepsilon}(x,s)-h_\varepsilon(x)}{\lambda})|+\sup_{s\in[T,t]}|u^{\lambda,N\varepsilon}(\pm N,s)-v(\pm N,s)|\right).
		\end{split}
\end{equation*} 
Letting $N$ go to $\infty$ and combining the uniform bound of $u^{\lambda,N,\varepsilon}$, it holds that 
$$
\sup_{x} |u^{\lambda,\varepsilon}(x,t)-v(x,t)| \le  C(\lambda,T)e^{r(T-t)}.
$$
Letting $t$ goes to $T$, we get that $u^{\lambda,\varepsilon}(x,0)=h_\varepsilon(x)$. Thus, $u^{\lambda,\varepsilon}$ is the solution of the following Cauchy problem
\begin{equation}\label{pde_epsilon}
	\left \{
	\begin{split}
		&\partial_t u^{\lambda,\varepsilon}+\mathcal L_x u^{\lambda,\varepsilon}+\lambda \exp(-\frac{u^{\lambda,\varepsilon}-h_\varepsilon}{\lambda})=0, (x,t)\in \mathbb R \times(0,T),\\
		&u^{\lambda,\varepsilon}(x,T)=h_\varepsilon(x).
	\end{split}
	\right.
\end{equation} 
and  satisfies the estimation \eqref{esti_bound}. 
}

For any compact set $D \subset \mathbb R \times(0,T)$, interior Schauder estimates yields that 
$$
\|u^{\lambda,\varepsilon}\|_{C^{2+\alpha,1+\frac{\alpha}{2}}(D)}\le C(\lambda,D)
$$
 Hence, $\{u^{\lambda,\varepsilon}\}$ locally uniformly converges to a function $u^\lambda$. Clearly, $u^\lambda$ satisfies the equation \eqref{lambda_HJB_2} in $\mathbb R\times(0,T)$. Let us check that it also satisfies the terminal condition. Denote by $\varphi(x,t)$ the value of European put option at time $t$ in our model setup with stock price $e^x$, which can be explicitly given, see \cite[Section 2.3]{jeanblanc2009mathematical}. It is also the solution for the following PDE
\begin{equation*}
	\left\{
	\begin{split}
	&\partial_t \varphi	+\mathcal L_x \varphi=0, (x,t) \in \mathbb R \times (0,T),\\
	&\varphi(x,T)=h(x).
	\end{split}
	\right.
\end{equation*}
Since $u^{\lambda,\varepsilon}$ is uniformly bounded, the nonlinear term $\lambda \exp(-\frac{u^{\lambda,\varepsilon}-h_\varepsilon}{\lambda})$ in \eqref{pde_epsilon} has a bound independent of $\varepsilon$. Using maximum principle again, it is not hard to prove that $|u^{\lambda,\varepsilon}(x,t)-\varphi(x,t)| \le C(\lambda)((T-t)+\|h-h_\varepsilon\|_\infty)$. Letting $\varepsilon$ go to zero and $t$ go to $T$, we verify the terminal condition for $u^\lambda$.}

\Copy{MathRigor2pr}{
\change{ Now, let us check that $u^\lambda(\log S,t)p$ is the value function $\tilde V(S,p,t)$. For any $\pi \in \mathcal A$, it holds that 
\begin{equation*}
	\small
	\begin{split}
	&e^{-rt}u^\lambda(\log S_t,t)p_t\\
	=&e^{-rT}h(log S_T)p_T-\int_t^Te^{-rs}\left(\frac12\sigma^2 \partial_{xx} u^\lambda+(r-\frac{1}{2}\sigma^2)\partial_x u-\pi_t u^\lambda-ru^\lambda\right)p_sds \\
	&-\int_t^T e^{-rs}\sigma p_s\partial_x u^\lambda dW_s\\
	=&e^{-rT}h(log S_T)p_T+\int_t^Te^{-rs}\left(\pi_t u^\lambda+ \lambda \exp(-\frac{u^\lambda-h}{\lambda})\right)p_s	ds -\int_t^T e^{-rs}\sigma p_s \partial_x u^\lambda dW_s.
\end{split}
\end{equation*}
Note that $\lambda\exp(-\frac{u^\lambda-h}{\lambda}) \ge (h-u^\lambda)\pi+\lambda R(\pi)$ for any $\pi>0$ and the equality holds only at $\pi=\exp(-\frac{u^\lambda-h}{\lambda})$.  Thus, 
\begin{equation*}
	\begin{split}
		&e^{-rt}u^\lambda(\log S_t,t)p_t\\
		\ge &e^{-rT}h(log S_T)p_T+\int_t^Te^{-rs}\left(h(log S_s)\pi_s+\lambda R(\pi_s) \right)p_s	ds -\int_t^T e^{-rs}\sigma p_s \partial_x u^\lambda dW_s.
	\end{split}
\end{equation*}
Taking conditional expectation, one shall get that $u^\lambda(\log S_t,t)p_t \ge J(p_t,S_t,t;\pi)$. For $\bar \pi^\lambda_t=\exp(-\frac{u^\lambda(\log S_t)-h(\log S_t)}{\lambda})$, all the inequality in the above will be equality and it holds that $u^\lambda(\log S_t,t)p_t \ge J(p_t,S_t,t;\bar \pi^\lambda)$. This implies that $u^\lambda(\log S,t)p$ is the optimal value function and $\pi^*$ is the optimal control.
}}
\subsection{Proofs for Theorem \ref{thm_monotone}}\Copy{MathProof4}{
Recall that $\{u^n,\pi^n\}$ are defined iteratively by the following equation
\begin{equation} \label{eq_VI_1}
	\partial_t u^n+\mathcal L_x u^n+H(x,\pi^n(x,t),u^n)=0,u^n(x,T)=h(x),
\end{equation}
and
\begin{equation}\label{eq_Policy_Update_1}
	\pi^{n+1}(x,t)=\mathop{\text{argmax}}_{\pi \in(0,\infty)}H(x,\pi,u^n(x,t))=\exp(-\frac{u^n(x,t)-h(x)}{\lambda}).
\end{equation}
\change{ Given $u^0 \in C_b(\mathbb R)$, we see that $\log \pi^1$ is also bounded as $\log \pi^1= -\frac{u^0-h}{\lambda}$ due to \eqref{eq_Policy_Update_1}. Then, from \eqref{eq_VI_1}, $u^1$ solves 
	$$
	\partial_t u^1+\mathcal L_x u^1+(h-u^1) \pi^1+\lambda R(\pi^1)=0,u^1(x,T)=h(x),
$$
which is a linear PDE with bounded coefficients. Using maximum principle, we obtain the boundedness of $u^1$. Repeating this argument, it is easy to obtain that $u^n$ is a bounded function for each $n$. Then, $\log \pi^n$ is also bounded for each $n$ as it is defined by \eqref{eq_Policy_Update_1}. 
}}
Define $\Delta^n:=u^{n+1}-u^{n}$. We see that it satisfies
\begin{equation}
	\partial_t \Delta^{n}+\mathcal L_x  \Delta^{n}+H(x,u^{n+1},\pi^{n+1})-H(x,u^n,\pi^n) =0,\Delta^{n}(x,T)=0.
\end{equation} 
From the definition of the Hamiltonian $H$, it can be also rewritten as 
\begin{equation}\label{pde_difference}
	\partial_t \Delta^{n}+\mathcal L_x  \Delta^{n}-\pi^{n+1}\Delta^n  =H(u^n,\pi^n)-H(u^n,\pi^{n+1})
\end{equation}
By the definition of $\pi^{n+1}$, we see that $H(u^n,\pi^n)-H(u^n,\pi^{n+1}) \le 0$. Hence, $0$ is a sub-solution  of \eqref{pde_difference}, which implies that $\Delta^n \ge 0$. Thus, we show that $u^{n+1} \ge u^{n}$. 

Next, consider the difference between $u^{\lambda}$ and $u^n$. Define $w^n:=u^\lambda-u^n$. Note that $u^\lambda$ satisfies
\begin{equation}\label{pde_lambda}
	\partial_t u^{\lambda}+\mathcal L_x  u^\lambda+H(x,u^\lambda,\pi^\lambda)=0,u^{\lambda}(x,T)=h(x).
\end{equation}
with 
$$
	\pi^{\lambda}(x,t)=\mathop{\text{argmax}}_{\pi \in(0,\infty)}H(x,\pi,u^\lambda(x,t))=\exp(-\frac{u^\lambda(x,t)-h(x)}{\lambda}).
$$
Thus, 
$$
\partial_t w^{n}+\mathcal L_x  w^n+H(x,u^\lambda,\pi^{\lambda})-H(x,u^{n},\pi^n)=0,w^n(x,T)=0.
$$
This is equivalent to  
\begin{equation}\label{eq_diff}
	\partial_t w^{n}+\mathcal L_x  w^n- \pi^n w^n =H(u^\lambda,\pi^n)-H(u^\lambda,\pi^\lambda)
\end{equation}
By the definition of $\pi^\lambda$, the right hand side is less than $0$. Applying maximum principle, we shall have that $w^n \ge 0$. At this stage, we have proved that $u^{n} \le u^{\lambda}$.
\subsection{Proof for Theorem \ref{thm_convergence_rate}} \Copy{MathProof7}{
 We have shown that $u^1 \le u^n \le u^{\lambda}$. Hence, $\{ u^n\}$ is uniformly bounded. 
 Moreover, $u^\lambda(x,t),u^n(x,t) \ge h(x)-M$ with $M:=\|(h-u^1)^+\|_\infty$. To prove the convergence rate, we define a function $f^n(t)$ as 
$$
f^n(t):=\sup_{x} u^{\lambda}(x,t)-u^n(x,t).
$$ 
By the definition of $\pi^\lambda$ and $\pi^n$, we compute that 
\begin{equation*}
	\begin{split}
		&H(x,u^\lambda,\pi^{n+1})-H(x,u^{\lambda},\pi^\lambda)\\
		=& (u^n-u^\lambda) \exp(-\frac{u^n-h}{\lambda})+\lambda(\exp(-\frac{u^n-h}{\lambda})-\exp(-\frac{u^\lambda-h}{\lambda}) ) 
	\end{split}
\end{equation*}
\change{
Then, it holds that 
$$
|H(x,u^\lambda,\pi^{n+1})-H(x,u^{\lambda},\pi^\lambda)| \le C^*\exp(\frac{M}{\lambda}) |u^n(x,t)-u^\lambda(x,t)|\le  C^*\exp(\frac{M}{\lambda}) f^n(t),
$$
with a constant $C^*$ independent of $\lambda$ and $n$.} Define $\phi(t):=C^*\exp(\frac{M}{\lambda}+rt)\int_t^T f^{n}(s)ds$. It solves $\partial_t \phi-r\phi=-C^*\exp(\frac{M}{\lambda}+rt)f^{n}(t)$. Then, we see that $W:=u^\lambda-u^{n+1}-\phi$ satisfies
$$
\partial_t W+\mathcal L_x  W=(u^\lambda-u^{n+1}) \pi^{n+1}+H(u^\lambda,\pi^{n+1})-H(u^\lambda,\pi^\lambda)+C^*\exp(\frac{M}{\lambda}+rt)f^{n}(t),
$$
with terminal condition $W(x,T)=0$. Previous argument shows that $$H(u^\lambda,\pi^{n+1})-H(u^\lambda,\pi^\lambda)+C^*\exp(\frac{M}{\lambda}+rt)f^{n}(t) \ge 0 .$$ Combining with the fact that $u^\lambda-u^{n+1} \ge 0$, it holds that 
$$
\partial_t W+\mathcal L_x  W \ge 0,
$$ 
which will implies that $W \le 0$ from maximum principle. This also means that 
$$
u^{\lambda}(x,t)-u^{n+1}(x,t) \le \phi(t)=C^*\exp(\frac{M}{\lambda}+rt) \int_t^T f^{n}(s)ds.
$$
Taking supreme with respect to $x$ on the left hand side, we have that 
$$
f^{n+1}(t) \le C^*\exp(\frac{M}{\lambda}+rt) \int_t^T f^{n}(s)ds \le C^*e^{rT}\exp(\frac{M}{\lambda}) \int_t^T f^{n}(s)ds .
$$
It is easy to see that $f^1$ is bounded in $t$. Proving by induction, we see that such that $f^{n+1}(t) \le (C^*)^n\frac{(T-t)^n}{n!}\exp(n(\frac{M}{\lambda}+rT))\sup_t|f^1(t)|$. Finally, we see that 
$$
|f^1(t)| =\sup_x u^\lambda(x,t)-u^1(x,t) \le \sup_x  u^\lambda(x,t)-h(x,t) +\sup_x (h-u^1)^+(x,t).
$$
Hence, $\sup_t |f^1(t)|\le \|u^\lambda-h\|_\infty+\|(u^1-h)^+\|_{\infty}\le M+K+\lambda T$. Hence, we finish the proof.}
\subsection{Proof for Lemma \ref{lem_lower_bound}}
Let $\kappa$ be a smooth, increasing, convex function such that $\kappa(x)=x$ for  $x\ge 1$ and $\kappa(x)=0$ for $x\le 0$. Denote by $C_1$, the supreme of $\kappa'$.  Given $\varepsilon>0$, define $\kappa_\varepsilon(x):=\varepsilon\kappa(\frac{K-e^x}{\varepsilon})$. It is easy to check that $\|\kappa_\varepsilon(x)-h(x)\|_\infty \rightarrow 0$ as $\varepsilon \rightarrow 0$. 

For any $\delta>0$, let us prove that $u^\lambda(x,t) \ge h(x)-\lambda \log \frac{(2+C_1) rK+\delta}{\lambda}$. Assume the contrary, i.e. $\inf_{(x,t)\in \mathbb R\times [0,T]} u^\lambda(x,t)-(h(x)-\lambda \log\frac{(2+C_1)rK+\delta}{\lambda})<0$. Since $\kappa_\varepsilon$ uniformly convergent to $h(x)$, we also have that $\inf_{(x,t)\in \mathbb R\times [0,T]} u^\lambda(x,t)-(\kappa_{\varepsilon_1}(x)-\lambda \log\frac{(2+C_1)rK+\delta}{\lambda})<0$ for sufficiently small ${\varepsilon_1}$. Choose another sufficient small constant $\varepsilon_2$. It will hold that 
$ \inf_{(x,t)\in \mathbb R\times [0,T]} u^\lambda(x,t)-(\kappa_{\varepsilon_1}(x)-\lambda \log\frac{(2+C_1)rK+\delta}{\lambda}-\frac{\varepsilon_2}{2} e^{rt} (|x|^2+\frac{1}{t}))<0$ and the infimum is achieved at some point $(x^*,t^*)$ with $t^*\in (0,T)$. \change{Then, it holds that 
$$
\frac{\varepsilon_2}{2} e^{rt^*} (|x^*|^2+\frac{1}{t^*})) \le \kappa_{\varepsilon_1}(x^*)-u^\lambda(x^*,t^*)-\lambda \log\frac{(2+C_1)rK+\delta}{\lambda}.
$$
Note that the third term $\lambda \log\frac{(2+C_1)rK+\delta}{\lambda}$ goes to $0$ as $\lambda$ goes to $0$ and goes to $\infty$ as $\lambda$ goes to $\infty$ and hence lower bounded with respect to $\lambda$.
} Since $u^\lambda$ and $\kappa_{\varepsilon_1}$ are bounded, we shall have that $\frac{\varepsilon_2}{2}e^{rt^*}(|x^*|^2+\frac{1}{t^*}) \le C+K+\lambda T$ with the constant $C$ independent of $\varepsilon_1$ and $\varepsilon_2$. Denote by $w(x,t)=\kappa_{\varepsilon_1}(x)-\lambda \log\frac{(2+C_1)rK+\delta}{\lambda}-\frac{\varepsilon_2}{2} e^{rt}(|x|^2+\frac{1}{t})$. It holds that 
$$
\partial_x w= -\kappa'(\frac{K-e^x}{\varepsilon_1})e^x-\varepsilon_2 e^{rt} x,
$$
$$
\partial_{xx}w=
\frac{1}{\varepsilon_1} \kappa''(\frac{K-e^x}{\varepsilon_1})e^{2x}-\kappa'(\frac{k-e^x}{\varepsilon_1})e^x-\varepsilon_2 e^{rt},
$$
and $\partial_t w=\frac{\varepsilon_2}{2}\frac{e^{rt}}{t^2}-r\frac{\varepsilon_2}{2}\frac{e^{rt}}{t} $. At $(x^*,t^*)$, we have 
$$
\partial_x u^\lambda(x^*,t^*)=\partial_{x} w(x^*,t^*),\partial_{xx} u^\lambda(x^*,t^*)\ge\partial_{xx} w(x^*,t^*),$$
$ \partial_t u^\lambda(x^*,t^*)=\partial_{t} w(x^*,t^*)$ and $u^\lambda(x^*,t^*) < w(x^*,t^*)$.  Then, at $(x^*,t^*)$
 \begin{equation}\label{eq_3.6.1}
 	\begin{split}
 	0=&\partial_t u^\lambda+\mathcal L_x u^\lambda+\lambda \exp(-\frac{u^\lambda(x^*,t^*) -h(x^*)}{\lambda})\\
 	 >& \partial_t w+\mathcal L_x w+\lambda \exp(-\frac{u^\lambda(x^*,t^*) -h(x^*)}{\lambda}).
 	\end{split}
 \end{equation}
Since $\kappa$ is convex, we have 
 $$
 \partial_t w(x^*,t^*)+\mathcal L_x w(x^*,t^*) \ge -r\kappa'(\frac{K-e^{x^*}}{\varepsilon_1})e^{x^*} -\varepsilon_2(r-\frac{\sigma^2}{2}) e^{rt^*}x^*-\frac{\sigma^2}{2}\varepsilon_2 e^{rt^*}-r\kappa_{\varepsilon_1}(x^*).
 $$
 Note that $\kappa'(\frac{K-e^{x^*}}{\varepsilon_1})$ is non-zero only if $e^{x^*} \le K$ and $\frac{\varepsilon_2}{2}e^{rt^*} |x^*|^2 \le C+K+\lambda T$. Since $\kappa_{\varepsilon_1}(x) \le K$, we shall have that 
 \begin{equation}\label{eq_3.6.2}
 \partial_t w(x^*,t^*)+\mathcal L_x w(x^*,t^*) \ge -(1+C_1)rK -\sqrt{\varepsilon_2}C(\lambda,r,T,K).
 \end{equation}
At $(x^*,t^*)$, it holds that 
 \begin{equation*}
 	\begin{split}
 		&h(x^*)-u^\lambda(x^*,t^*)\\
 		\ge& h(x^*)-w(x^*,t^*) \\
 		\ge& h(x^*)-\kappa_{\varepsilon_1}(x^*)+\lambda\log\frac{(2+C_1)rK+\delta}{\lambda}.
 		\end{split}
 	\end{equation*}
 Thus, 
 \begin{equation}\label{eq_3.6.3}
 	 \lambda \exp(-\frac{u^\lambda(x^*,t^*)- h(x^*)}{\lambda})\ge ((2+C_1)rK+\delta)\exp(\frac{h(x^*)-\kappa_{\varepsilon_1}(x^*)}{\lambda})
 \end{equation}
Combining \eqref{eq_3.6.2} and  \eqref{eq_3.6.3}, we see that the right hand side of \eqref{eq_3.6.1} will be strictly positive for $\varepsilon_1$ and $\varepsilon_2$ sufficiently small. This leads to a contradiction. Thus, we prove that, for any $\delta>0$, $u^\lambda(x,t) \ge h(x)-\lambda \log \frac{(2+C_1)rK+\delta}{\lambda}$. Letting $\delta \rightarrow 0$, we have the desired result.
\subsection{Proof for Theorem \ref{thm_difference}}\Copy{MathProof5}{
\change{ In \cite[Theorem 3.6]{jaillet1990variational}, the authors proved the existence of $u^*$ and showed that $\partial_t u,\partial_x u$ and $\partial_{xx} u$ are locally bounded in $\mathbb R \times [0,T)$. In the region $\{ u^*>h\}$, $u^*$ satisfies $\partial_t u^* +\mathcal L_x u^*=0$. Hence, applying local Schauder estimate, one can show that the derivatives are also continuous.
}}

\Copy{MathProof6}{ Assume the contrary that $\sup_{(x,t) \in \mathbb R \times[0,T]} u^*(x,t)-u^{\lambda}(x,t)-\lambda \log\frac{CrK}{\lambda} >0$. \change{ This means that there exists $(x',t') \in \mathbb R \times (0,T)$ such that $u^*(x',t')-u^{\lambda}(x',t')-\lambda \log\frac{CrK}{\lambda}>0$. Then, choosing $\varepsilon $ sufficiently small, it holds that $u^*(x',t')-u^{\lambda}(x',t')-\lambda \log\frac{CrK}{\lambda}-\varepsilon (\frac12|x'|^2+\frac1t')>0$. Thus, one can choose $\varepsilon$ sufficiently small such that $\sup_{(x,t) \in \mathbb R \times[0,T]} u^*(x,t)-u^{\lambda}(x,t)-\lambda \log\frac{CrK}{\lambda}-\varepsilon(\frac12|x|^2+\frac1t)>0$.  Since $\frac12|x|^2+\frac1t$ goes to $\infty$ when $|x| \rightarrow \infty$ and $t \rightarrow 0$, the supreme will be attained at some point $(x^*,t^*) \in \mathbb R \times(0,T)$.} }
 By the terminal condition of $u^*$ and $u^\lambda$, it holds that $t^* \in (0,T)$. Denote by $w=u^\lambda +\lambda \log \frac{CrK}{\lambda}+\varepsilon(\frac12|x|^2+\frac1t)$. Then, at $(x^*,t^*)$, it holds that 
$$
\partial_{xx} u^* (x^*,t^*) \le \partial_{xx} w (x^*,t^*),\partial_{x} u^* (x^*,t^*) = \partial_{x} w (x^*,t^*) \text{ and } \partial_{t} u^* (x^*,t^*) = \partial_{t} w.
$$
We also have that 
$$
 \partial_{xx} w (x^*,t^*)= \partial_{xx} u^\lambda (x^*,t^*)+\varepsilon, \partial_{x} w (x^*,t^*)= \partial_{x} u^\lambda (x^*,t^*)+\varepsilon x^*,
$$
and $\partial_{t} w (x^*,t^*)= \partial_{t} u^\lambda (x^*,t^*)-\varepsilon \frac{1}{(t^*)^2}$. Since  $u^*(x^*,t^*) > u^{\lambda}(x^*,t^*)+\lambda \log\frac{CrK}{\lambda}$ and Lemma 1, we have $u^*(x^*,t^*)>h(x^*)$. This implies that, at $(x^*,t^*)$ 
$$
\partial_t u^* +\mathcal L_x u^*=0. 
$$ 
Then, from the relation between the derivatives of $u^*$ and $w$, we deduce that, at $(x^*,t^*)$
$$
\partial_t u^\lambda +\mathcal L_x u^\lambda\ge -\frac{1}{2}\sigma^2 \varepsilon -(r-\frac{1}{2}\sigma^2)\varepsilon x^*+\varepsilon \frac{1}{(t^*)^2}.
$$
We also have $\frac{\varepsilon}{2}|x^*|^2 \le \tilde C$ with the constant $\tilde C$ independent of $\varepsilon$. This implies that $\varepsilon x^*=O(\varepsilon^{1/2})$. In Theorem \ref{thm_solvability}, we have shown that $u^\lambda$ is bounded, which yields $\lambda \exp(-\frac{u^\lambda-h}{\lambda})$ is uniformly larger than some small constant.  Hence, we shall have that  $\partial_t u^\lambda +\mathcal L_x u^\lambda+\lambda \exp(-\frac{u^{\lambda}-h}{\lambda})>0$ at $(x^*,t^*)$ for $\varepsilon$ sufficiently small. This leads to a contradiction. Hence, we have prove the second inequality.

Let $f(t):=\log(1+(T-t))$. It is easy to see that 
$$
\partial_t f+\exp(-f)=0.
$$
Denote by $w(x,t)=u^*(x,t)+\lambda f(t)$. Then, we see that $\exp(-\frac{w-h}{\lambda}) \le \exp(-f)$ as $u^*\ge h$. Combining with the fact that $\partial_t u^*+\mathcal L_x u^* \le 0$, it holds that $w$ solves the following inequality
$$
\partial_t w+\mathcal L_x w +\lambda \exp(-\frac{w-h}{\lambda}) < 0,
$$
with terminal condition $w(x,T)=u^\lambda(x,T)$. Then, applying maximum principle, we shall have that 
$u^\lambda \le w$, which is equivalent to the first inequality of our theorem. 
\subsection{Proof for Theorem \ref{thm_wealth_loss}}
Applying It\^o formula, we see that 
\begin{equation*}
	\begin{split}
		&e^{-r\tau^\lambda}g(S_{\tau^{\lambda}})\\
		=&e^{-r\tau^\lambda}V^{\lambda}(S_{\tau^{\lambda}},\tau^{\lambda})\\
		=&V^{\lambda}(S_0,0)+\int_0^{\tau^{\lambda}} e^{-ru}(\partial_t+\mathcal L_S )V^{\lambda}(S_u,u)du+\sigma\int_0^{\tau^{\lambda}} e^{-ru} \partial_S V(S_u,u)dW_u\\
		=&V^{\lambda}(S_0,0)-\lambda\int_0^{\tau^{\lambda}}e^{-ru}\exp(-\frac{V^{\lambda}(S_u,u)-g(S_u)}{\lambda})du+\sigma\int_0^{\tau^{\lambda}} e^{-ru}\partial_S V(S_u,u)dW_u\\
		\ge&V^{\lambda}(S_0,0)-\lambda T+\sigma\int_0^{\tau^{\lambda}} e^{-ru}\partial_S V(S_u,u)dW_u.
	\end{split}
\end{equation*}
Taking expectation, we have 
$$
\mathbb E\left[ e^{-r\tau^\lambda}	g(S_{\tau^{\lambda}}) \right] \ge V^{\lambda}(S_0,0)-\lambda T \ge V^*(S_0,0)+ \inf_S (V^{\lambda}(\cdot,0)-V^*(\cdot,0))-\lambda T.
$$
Previous lemma indicates that  $ \inf_S (V^{\lambda}(\cdot,0)-V^*(\cdot,0)) \ge -\lambda \log{\frac{CrK}{\lambda}}$. Hence, 
$$
\mathbb E\left[ e^{-r\tau^\lambda}	g(S_{\tau^{\lambda}}) \right] \ge V^{*}(S_0,0)-(\lambda\log\frac{CrK}{\lambda}+\lambda T).
$$
This implies that using learned strategy leads to a loss of wealth approximately of the order $O(\lambda \log\frac{1}{\lambda})$.
\bibliographystyle{siamplain}
\bibliography{references}

\begin{thebibliography}{10}

\bibitem{ajdari2019towards}
{\sc A.~Ajdari, M.~Niyazi, N.~H. Nicolay, C.~Thieke, R.~Jeraj, and
  T.~Bortfeld}, {\em Towards optimal stopping in radiation therapy},
  Radiotherapy and Oncology, 134 (2019), pp.~96--100.

\bibitem{becker2019deep}
{\sc S.~Becker, P.~Cheridito, and A.~Jentzen}, {\em Deep optimal stopping},
  Journal of Machine Learning Research, 20 (2019), p.~74.

\bibitem{becker2021solving}
{\sc S.~Becker, P.~Cheridito, A.~Jentzen, and T.~Welti}, {\em Solving
  high-dimensional optimal stopping problems using deep learning}, European
  Journal of Applied Mathematics, 32 (2021), pp.~470--514.

\bibitem{dai2019bayesian}
{\sc Z.~Dai, H.~Yu, B.~K.~H. Low, and P.~Jaillet}, {\em Bayesian optimization
  meets bayesian optimal stopping}, in International conference on machine
  learning, PMLR, 2019, pp.~1496--1506.

\bibitem{deisenroth2013survey}
{\sc M.~P. Deisenroth, G.~Neumann, J.~Peters, et~al.}, {\em A survey on policy
  search for robotics}, Foundations and Trends{\textregistered} in Robotics, 2
  (2013), pp.~1--142.

\bibitem{han2017deep}
{\sc W.~E, J.~Han, and A.~Jentzen}, {\em Deep learning-based numerical methods
  for high-dimensional parabolic partial differential equations and backward
  stochastic differential equations}, Communications in mathematics and
  statistics, 5 (2017), pp.~349--380.

\bibitem{fathan2021deep}
{\sc A.~Fathan and E.~Delage}, {\em Deep reinforcement learning for optimal
  stopping with application in financial engineering}, arXiv preprint
  arXiv:2105.08877,  (2021).

\bibitem{felizardo2022solving}
{\sc L.~K. Felizardo, E.~Matsumoto, and E.~Del-Moral-Hernandez}, {\em Solving
  the optimal stopping problem with reinforcement learning: an application in
  financial option exercise}, in 2022 International Joint Conference on Neural
  Networks (IJCNN), IEEE, 2022, pp.~1--8.

\bibitem{forsyth2002quadratic}
{\sc P.~A. Forsyth and K.~R. Vetzal}, {\em Quadratic convergence for valuing
  american options using a penalty method}, SIAM Journal on Scientific
  Computing, 23 (2002), pp.~2095--2122.

\bibitem{fox2015taming}
{\sc R.~Fox, A.~Pakman, and N.~Tishby}, {\em Taming the noise in reinforcement
  learning via soft updates}, arXiv preprint arXiv:1512.08562,  (2015).

\bibitem{friedman2008partial}
{\sc A.~Friedman}, {\em Partial differential equations of parabolic type},
  Courier Dover Publications, 2008.

\bibitem{haarnoja2017reinforcement}
{\sc T.~Haarnoja, H.~Tang, P.~Abbeel, and S.~Levine}, {\em Reinforcement
  learning with deep energy-based policies}, in International conference on
  machine learning, PMLR, 2017, pp.~1352--1361.

\bibitem{hendricks2014reinforcement}
{\sc D.~Hendricks and D.~Wilcox}, {\em A reinforcement learning extension to
  the {A}lmgren-{C}hriss framework for optimal trade execution}, in 2014 IEEE
  Conference on Computational Intelligence for Financial Engineering \&
  Economics (CIFEr), IEEE, 2014, pp.~457--464.

\bibitem{hu2011blow}
{\sc B.~Hu}, {\em Blow-up theories for semilinear parabolic equations},
  Springer, 2011.

\bibitem{jaillet1990variational}
{\sc P.~Jaillet, D.~Lamberton, and B.~Lapeyre}, {\em Variational inequalities
  and the pricing of american options}, Acta Applicandae Mathematicae, 21
  (1990), pp.~263--289.

\bibitem{jeanblanc2009mathematical}
{\sc M.~Jeanblanc, M.~Yor, and M.~Chesney}, {\em Mathematical methods for
  financial markets}, Springer Science \& Business Media, 2009.

\bibitem{karniadakis2021physics}
{\sc G.~E. Karniadakis, I.~G. Kevrekidis, L.~Lu, P.~Perdikaris, S.~Wang, and
  L.~Yang}, {\em Physics-informed machine learning}, Nature Reviews Physics, 3
  (2021), pp.~422--440.

\bibitem{kerimkulov2020exponential}
{\sc B.~Kerimkulov, D.~Siska, and L.~Szpruch}, {\em Exponential convergence and
  stability of howard's policy improvement algorithm for controlled
  diffusions}, SIAM Journal on Control and Optimization, 58 (2020),
  pp.~1314--1340.

\bibitem{kulis2009low}
{\sc B.~Kulis, M.~A. Sustik, and I.~S. Dhillon}, {\em Low-rank kernel learning
  with {B}regman matrix divergences.}, Journal of Machine Learning Research, 10
  (2009).

\bibitem{laraki2005value}
{\sc R.~Laraki and E.~Solan}, {\em The value of zero-sum stopping games in
  continuous time}, SIAM Journal on Control and Optimization, 43 (2005),
  pp.~1913--1922.

\bibitem{levine2016end}
{\sc S.~Levine, C.~Finn, T.~Darrell, and P.~Abbeel}, {\em End-to-end training
  of deep visuomotor policies}, The Journal of Machine Learning Research, 17
  (2016), pp.~1334--1373.

\bibitem{lieberman1996second}
{\sc G.~M. Lieberman}, {\em Second order parabolic differential equations},
  World scientific, 1996.

\bibitem{liyanage2019automating}
{\sc Y.~W. Liyanage, D.-S. Zois, C.~Chelmis, and M.~Yao}, {\em Automating the
  classification of urban issue reports: an optimal stopping approach}, in
  ICASSP 2019-2019 IEEE International Conference on Acoustics, Speech and
  Signal Processing (ICASSP), IEEE, 2019, pp.~3137--3141.

\bibitem{longstaff2001valuing}
{\sc F.~A. Longstaff and E.~S. Schwartz}, {\em Valuing american options by
  simulation: a simple least-squares approach}, The review of financial
  studies, 14 (2001), pp.~113--147.

\bibitem{mnih2015human}
{\sc V.~Mnih, K.~Kavukcuoglu, D.~Silver, A.~A. Rusu, J.~Veness, M.~G.
  Bellemare, A.~Graves, M.~Riedmiller, A.~K. Fidjeland, G.~Ostrovski, et~al.},
  {\em Human-level control through deep reinforcement learning}, Nature, 518
  (2015), p.~529.

\bibitem{moody2001learning}
{\sc J.~Moody and M.~Saffell}, {\em Learning to trade via direct
  reinforcement}, IEEE transactions on neural Networks, 12 (2001),
  pp.~875--889.

\bibitem{moody1998performance}
{\sc J.~Moody, L.~Wu, Y.~Liao, and M.~Saffell}, {\em Performance functions and
  reinforcement learning for trading systems and portfolios}, Journal of
  Forecasting, 17 (1998), pp.~441--470.

\bibitem{nachum2017bridging}
{\sc O.~Nachum, M.~Norouzi, K.~Xu, and D.~Schuurmans}, {\em Bridging the gap
  between value and policy based reinforcement learning}, in Advances in Neural
  Information Processing Systems, 2017, pp.~2775--2785.

\bibitem{nevmyvaka2006reinforcement}
{\sc Y.~Nevmyvaka, Y.~Feng, and M.~Kearns}, {\em Reinforcement learning for
  optimized trade execution}, in Proceedings of the 23rd international
  conference on Machine learning, ACM, 2006, pp.~673--680.

\bibitem{raissi2017physics}
{\sc M.~Raissi, P.~Perdikaris, and G.~E. Karniadakis}, {\em Physics informed
  deep learning (part i): Data-driven solutions of nonlinear partial
  differential equations}, arXiv preprint arXiv:1711.10561,  (2017).

\bibitem{reppen2022neural}
{\sc A.~M. Reppen, H.~M. Soner, and V.~Tissot-Daguette}, {\em Neural optimal
  stopping boundary}, arXiv preprint arXiv:2205.04595,  (2022).

\bibitem{silver2016mastering}
{\sc D.~Silver, A.~Huang, C.~J. Maddison, A.~Guez, L.~Sifre, G.~Van
  Den~Driessche, J.~Schrittwieser, I.~Antonoglou, V.~Panneershelvam,
  M.~Lanctot, et~al.}, {\em Mastering the game of {G}o with deep neural
  networks and tree search}, nature, 529 (2016), p.~484.

\bibitem{silver2017mastering}
{\sc D.~Silver, J.~Schrittwieser, K.~Simonyan, I.~Antonoglou, A.~Huang,
  A.~Guez, T.~Hubert, L.~Baker, M.~Lai, A.~Bolton, et~al.}, {\em Mastering the
  game of go without human knowledge}, Nature, 550 (2017), p.~354.

\bibitem{sirignano2018dgm}
{\sc J.~Sirignano and K.~Spiliopoulos}, {\em {DGM}: A deep learning algorithm
  for solving partial differential equations}, Journal of computational
  physics, 375 (2018), pp.~1339--1364.

\bibitem{smirnov}
{\sc M.~Smirnov}, {\em Javascript options and implied volatility calculator.}
\newblock \url{http://www.math.columbia.edu/~smirnov/options13.html}.

\bibitem{sutton2011reinforcement}
{\sc R.~S. Sutton and A.~G. Barto}, {\em Reinforcement learning: An
  introduction}, Cambridge, MA: MIT Press, 2011.

\bibitem{tang2021exploratory}
{\sc W.~Tang, P.~Y. Zhang, and X.~Y. Zhou}, {\em Exploratory {HJB} equations
  and their convergence}, arXiv preprint arXiv:2109.10269,  (2021).

\bibitem{touzi2002continuous}
{\sc N.~Touzi and N.~Vieille}, {\em Continuous-time dynkin games with mixed
  strategies}, SIAM Journal on Control and Optimization, 41 (2002),
  pp.~1073--1088.

\bibitem{wang2018exploration}
{\sc H.~Wang, T.~Zariphopoulou, and X.~Zhou}, {\em Exploration versus
  exploitation in reinforcement learning: a stochastic control approach}, arXiv
  preprint arXiv:1812.01552,  (2018).

\bibitem{wang2020reinforcement}
{\sc H.~Wang, T.~Zariphopoulou, and X.~Y. Zhou}, {\em Reinforcement learning in
  continuous time and space: A stochastic control approach.}, J. Mach. Learn.
  Res., 21 (2020), pp.~198--1.

\bibitem{wang2021deep}
{\sc S.~Wang and P.~Perdikaris}, {\em Deep learning of free boundary and
  {S}tefan problems}, Journal of Computational Physics, 428 (2021), p.~109914.

\bibitem{ziebart2008maximum}
{\sc B.~D. Ziebart, A.~Maas, J.~A. Bagnell, and A.~K. Dey}, {\em Maximum
  entropy inverse reinforcement learning},  (2008).

\bibitem{zimmert2019connections}
{\sc J.~Zimmert and T.~Lattimore}, {\em Connections between mirror descent,
  thompson sampling and the information ratio}, Advances in Neural Information
  Processing Systems, 32 (2019).

\end{thebibliography}

\end{document}